\documentclass{amsart}

\usepackage{amsmath}
\usepackage{amssymb}

\newtheorem{theorem}{Theorem}[section]
\newtheorem{lemma}[theorem]{Lemma}
\newtheorem{proposition}[theorem]{Proposition}
\newtheorem{corollary}[theorem]{Corollary}

\theoremstyle{definition}
\newtheorem{definition}[theorem]{Definition}

\theoremstyle{remark}
\newtheorem{remark}[theorem]{Remark}

\numberwithin{equation}{section}

\begin{document}
\title[Biharmonic homogeneous hypersurfaces]{Biharmonic homogeneous hypersurfaces in compact symmetric spaces}

\author{Shinji Ohno}
\address{
Department of Mathematics and Information Sciences, 
Tokyo Metropolitan University, 
Minami-Osawa 1-1, Hachioji, Tokyo, 192-0397, Japan.}
\email{oono-shinji@ed.tmu.ac.jp}

\author{Takashi Sakai}
\address{
Department of Mathematics and Information Sciences, 
Tokyo Metropolitan University, 
Minami-Osawa 1-1, Hachioji, Tokyo, 192-0397, Japan.}
\email{sakai-t@tmu.ac.jp}

\author{Hajime Urakawa}
\address{
Institute for International Education, 
Tohoku University, 
Kawauchi 41, Sendai, 980-8576, Japan.}
\email{urakawa@math.is.tohoku.ac.jp}

\thanks{
The second author was supported by Grant-in-Aid for Scientific Research (C) No.~26400073,
Japan Society for the Promotion of Science.
The third author was supported by Grant-in-Aid for Scientific Research (C) No.~25400154,
Japan Society for the Promotion of Science.
}


\date{\today}


\begin{abstract}
In this paper, we study biharmonic hypersurfaces in Einstein manifolds.
Then, we determine all the biharmonic hypersurfaces
in irreducible symmetric spaces of compact type
which are regular orbits of commutative Hermann actions of cohomogeneity one. 
\end{abstract}

\maketitle

\section{Introduction} \label{sect:introduction}
Harmonic maps play a central role in geometry;\,they are critical points of the energy functional 
$E(\varphi)=(1/2)\int_M\vert d\varphi\vert^2\,v_g$ 
for smooth maps $\varphi$ of $(M,g)$ into $(N,h)$. The Euler-Lagrange equations are given by the vanishing of the tension filed 
$\tau(\varphi)$. 
In 1983, J. Eells and L. Lemaire \cite{EL1} extended the notion of harmonic map to  
biharmonic map, which are, 
by definition, 
critical points of the bienergy functional
\begin{equation*}
E_2(\varphi)=\frac12\int_M
\vert\tau(\varphi)\vert^2\,v_g.
\end{equation*}
After G.Y. Jiang \cite{J} studied the first and second variation formulas of $E_2$, 
extensive studies in this area have been done
(for instance, see 
\cite{CMP}, \cite{IIU2}, \cite{IIU},  \cite{II}, 
\cite{LO2},  \cite{MO1}, \cite{OT2}, \cite{S1}, etc.).
Notice that harmonic maps are always biharmonic by definition.
One of the important main problems is to ask whether the converse is true. 
B.Y. Chen raised (\cite{C}) so called B.Y. Chen's conjecture and later,
R. Caddeo, S. Montaldo, P. Piu and C. Oniciuc raised (\cite{CMP}) the generalized B.Y. Chen's conjecture: 
\par
{\em Every biharmonic submanifold of the Euclidean space $\mathbb{R}^n$ must be harmonic (minimal).}
\par
{\em Every biharmonic submanifold of a Riemannian manifold of non-positive curvature must be harmonic (minimal).}
\vskip0.2cm\par
For the generalized Chen's conjecture, 
Ou and Tang gave (\cite{OT}, \cite{OT2}) a counter example in a Riemannian manifold of negative curvature. 
The Chen's conjecture was solved affirmatively in the case of surfaces in the three dimensional Euclidean space (\cite{C}), 
and the case of hypersurfaces of the four dimensional Euclidean space (\cite{D}, \cite{HV}),
and the case of generic hypersurfaces in the Euclidean space (\cite{KU}).
\par
Furthermore, Akutagawa and Maeta gave (\cite{AM}) a final supporting evidence to the Chen's conjecture: 
Every complete properly immersed biharmonic submanifold of the Euclidean space 
$\mathbb{R}^n$ is minimal. 
\par
It is also known (cf. \cite{NU1}, \cite{NU2}, \cite{NUG}): 
every biharmonic map $\varphi:\,(M,g)\rightarrow (N,h)$ of a complete Riemannian manifold $(M,g)$
into another Riemannian manifold $(N,h)$ with non-positive sectional curvature
with finite energy and finite bienergy is harmonic. 
\par
On the contrary to the above, the case that the target space $(N,h)$ whose sectional curvature is non-negative,
theory of biharmonic maps and/or biharmonic immersions is quite different. 
In 1986, Jiang \cite{J} and in 2002, Oniciuc \cite{O} constructed independently different examples of proper biharmonic
immersions into the spheres. 
Here, {\em proper biharmonic} means that biharmonic, but not harmonic. 
\par
In this paper, we study biharmonic hypersurfaces in the Einstein manifold (cf. Theorems~\ref{Theorem 3.6} and \ref{thm:3.7}),
and then we determine all the biharmonic hypersurfaces in the irreducible symmetric spaces of compact type
which are regular orbits of commutative Hermann actions of cohomogeneity one (cf. Theorem~~\ref{thm:list_of_biharmonic_orbits}).
More precisely, we obtain the following 
(cf. Theorem~~\ref{thm:list_of_biharmonic_orbits}).
\par
Let us consider all the commutative compact symmetric triads 
$(G,K_1,K_2)$ whose Hermann actions are of cohomogeneity one. 
Then, it holds that
an orbit $K_2\pi_1(x)$ in $N_1$ is harmonic (resp. proper biharmonic) if and only if the orbit 
$K_1\pi_2(x)$ in $N_2$ is harmonic (resp. proper biharmonic), 
where $\pi_i$ are the projections of $G$ onto $N_i=G/K_i$ $(i=1,2)$ (cf. Proposition 5.1).  
\par Furthermore, when $G$ is simple,
all the regular orbits of $K_2$ actions on $N_1=G/K_1$, and $K_1$ actions on $N_2=G/K_2$ are classified 
into the following three classes of totally eighteen 
cases of actions on symmetric space
(cf. Theorem~~\ref{thm:list_of_biharmonic_orbits}):
\begin{enumerate}
\item The first class consists of three cases.
Every case in this class has a unique regular orbit which is proper biharmonic. 
\item The second class consists of seven cases.
Each case in this class has exactly two regular orbits which are proper biharmonic. 
\item The third class consists of eight cases.
All the biharmonic regular orbits in the cases in this class must be harmonic. 
\end{enumerate}
\vskip0.6cm\par
\section{Preliminaries} \label{sect:preliminaries}
We first prepare the materials for the first and second variational formulas for the bienergy functional and biharmonic maps. 
Let us recall the definition of a harmonic map $\varphi:\,(M,g)\rightarrow (N,h)$,
of a compact Riemannian manifold $(M,g)$ into another Riemannian manifold $(N,h)$, 
which is an extremal of the {\em energy functional} defined by 
$$
E(\varphi)=\int_Me(\varphi)\,v_g, 
$$
where $e(\varphi):=(1/2)\vert d\varphi\vert^2$ is called the energy density of $\varphi$.
That is, for any variation $\{\varphi_t\}$ of $\varphi$ with $\varphi_0=\varphi$,
\begin{equation} \label{eq:2.1}
\frac{d}{dt}\bigg\vert_{t=0}E(\varphi_t)=-\int_Mh(\tau(\varphi),V)v_g=0,
\end{equation}
where $V\in \Gamma(\varphi^{-1}TN)$ is a variation vector field along $\varphi$ which is given by 
$V(x)=\frac{d}{dt}\big\vert_{t=0}\varphi_t(x)\in T_{\varphi(x)}N$, $(x\in M)$, 
and the {\em tension field} of $\varphi$ is given by 
$\tau(\varphi)=\sum_{i=1}^m B_\varphi(e_i,e_i)\in \Gamma(\varphi^{-1}TN)$, 
where 
$\{e_i\}_{i=1}^m$ is a locally defined orthonormal frame field on $(M,g)$, 
and $B_\varphi$ is the second fundamental form of $\varphi$ defined by 
\begin{align*}
B_\varphi(X,Y)
&=(\widetilde{\nabla}d\varphi)(X,Y) \\
&=(\widetilde{\nabla}_Xd\varphi)(Y) \\
&=\overline{\nabla}_X(d\varphi(Y))-d\varphi(\nabla_XY),
\end{align*}
for all vector fields $X, Y\in {\frak X}(M)$. 
Here, 
$\nabla$, and $\nabla^h$ are Levi-Civita connections on $TM$, $TN$ of $(M,g)$, $(N,h)$, respectively,
and $\overline{\nabla}$, and $\widetilde{\nabla}$ are the induced ones on $\varphi^{-1}TN$, and $T^{\ast}M\otimes \varphi^{-1}TN$, respectively.
By (\ref{eq:2.1}), $\varphi$ is {\em harmonic} if and only if $\tau(\varphi)=0$. 
\par
The second variation formula is given as follows.
Assume that $\varphi$ is harmonic. 
Then, 
\begin{equation*}
\frac{d^2}{dt^2}\bigg\vert_{t=0}E(\varphi_t)=\int_Mh(J(V),V)v_g, 
\end{equation*}
where $J$ is an elliptic differential operator, called the {\em Jacobi operator} acting on $\Gamma(\varphi^{-1}TN)$ given by 
\begin{equation} \label{eq:2.2}
J(V)=\overline{\Delta}V-\mathcal{R}(V),
\end{equation}
where 
$\overline{\Delta}V
=\overline{\nabla}^{\ast}\overline{\nabla}V
=-\sum_{i=1}^m\{\overline{\nabla}_{e_i}\overline{\nabla}_{e_i}V-\overline{\nabla}_{\nabla_{e_i}e_i}V\}$ 
is the {\em rough Laplacian} and 
$\mathcal{R}$ is a linear operator on $\Gamma(\varphi^{-1}TN)$ given by 
$\mathcal{R}(V)=\sum_{i=1}^mR^h(V,d\varphi(e_i))d\varphi(e_i)$,
and $R^h$ is the curvature tensor of $(N,h)$ given by 
$R^h(U,V)W=\nabla^h_U(\nabla^h_VW)-\nabla^h_V(\nabla^h_UW)-\nabla^h_{[U,V]}W$ for $U,V,W\in {\frak X}(N)$. 
\par
J. Eells and L. Lemaire \cite{EL1} proposed polyharmonic ($k$-harmonic) maps and 
Jiang \cite{J} studied the first and second variation formulas of biharmonic maps.
Let us consider the {\em bienergy functional} defined by 
\begin{equation*}
E_2(\varphi)=\frac12\int_M\vert\tau(\varphi)\vert ^2v_g, 
\end{equation*}
where 
$\vert V\vert^2=h(V,V)$, $V\in \Gamma(\varphi^{-1}TN)$.  
\par
The first variation formula of the bienergy functional is given by
\begin{equation*}
\frac{d}{dt}\bigg\vert_{t=0}E_2(\varphi_t)
=-\int_Mh(\tau_2(\varphi),V)v_g.
\end{equation*}
Here, 
\begin{equation*}
\tau_2(\varphi)
:=J(\tau(\varphi))=\overline{\Delta}(\tau(\varphi))-\mathcal{R}(\tau(\varphi)),
\end{equation*}
which is called the {\em bitension field} of $\varphi$, and 
$J$ is given in (\ref{eq:2.2}).  
\par
A smooth map $\varphi$ of $(M,g)$ into $(N,h)$ is said to be 
{\em biharmonic} if 
$\tau_2(\varphi)=0$. 
By definition, every harmonic map is biharmonic. 
We say, for an immersion $\varphi:\,(M,g)\rightarrow (N,h)$  to be {\em proper biharmonic} if 
it is biharmonic but not harmonic. 
\vskip0.6cm\par

\section{Biharmonic isometric immersions} \label{sect:biharmonic_maps}
{\bf 3.1.} 
In the first part of this section, we first show a characterization theorem 
for an isometric immersion $\varphi$ 
of an $m$ dimensional Riemannian manifold $(M,g)$ into an $n$ dimensional Riemannian manifold $(N,h)$ 
whose tension field $\tau(\varphi)$ satisfies that 
$\overline{\nabla}^{\perp}_X\tau(\varphi)=0$ $(X\in {\frak X}(M))$ to be biharmonic,
where $\overline{\nabla}^{\perp}$ is the normal connection on the normal bundle $T^\perp M$. 
Let us recall 
the following theorem due to \cite{J}: 
\begin{theorem} 
Let $\varphi:\,(M^m,g)\rightarrow (N^n,h)$ be an isometric immersion. 
Assume that $\overline{\nabla}^{\perp}_X\tau(\varphi)=0$ for all 
$X\in{\frak X}(M)$. Then, $\varphi$ is biharmonic if and only if 
the following holds:
\begin{align}
-\sum_{j,k=1}^m&h\big(\tau(\varphi),R^h(d\varphi(e_j),d\varphi(e_k))d\varphi(e_k)\big)\,d\varphi(e_j)\nonumber\\
&+\sum_{j,k=1}^mh(\tau(\varphi),B_{\varphi}(e_j,e_k))\,B_{\varphi}(e_j,e_k)\nonumber\\
&\quad -\sum_{j=1}^mR^h(\tau(\varphi),d\varphi(e_j))\,d\varphi(e_j)=0, \label{eq:3.1}
\end{align}
where $\{e_j\}_{j=1}^m$ is a locally defined orthonormal frame field 
on $(M,g)$. 
\end{theorem}
\medskip\par
Here, 
let us apply the following general curvature tensorial properties 
(\cite{KN}, Vol. I, Pages 198, and 201)
to the first term of the left hand side of (\ref{eq:3.1}): 
\begin{align*}
h\big(W_1,R^h(W_3,W_4)W_2\big)=h\big(W_3,&R^h(W_1,W_2)W_4\big), \\
& (W_i\in{\frak X}(N),i=1,2,3,4).
\end{align*}
Then, we have 
\begin{align*}
h\big(\tau(\varphi),&R^h(d\varphi(e_j),d\varphi(e_k))d\varphi(e_k)\big) \\
&=h\big(d\varphi(e_j),R^h(\tau(\varphi),d\varphi(e_k))d\varphi(e_k)\big).
\end{align*}
Therefore, for the first term of (\ref{eq:3.1}), we have that 
\begin{align*}
\sum_{j=1}^mh\big(
d\varphi(e_j),\sum_{k=1}^mR^h(\tau(\varphi),d\varphi(e_k))d\varphi(e_k)\big)\,d\varphi(e_j)
\end{align*}
is equal to the tangential part of 
$\sum_{k=1}^mR^h(\tau(\varphi),d\varphi(e_k))\,d\varphi(e_k)$. 
Thus, the equation (\ref{eq:3.1}) is equivalent to
\begin{align}
-&\bigg(\sum_{k=1}^mR^h(\tau(\varphi),d\varphi(e_k))d\varphi(e_k)\bigg)^{\top}\nonumber\\
&+\sum_{j,k=1}^mh(\tau(\varphi),B_{\varphi}(e_j,e_k))\,B_{\varphi}(e_j,e_k)\nonumber\\
&\quad -\sum_{k=1}^mR^h(\tau(\varphi),d\varphi(e_k))\,d\varphi(e_k)=0, \label{eq:3.2}
\end{align}
where $W^{\top}$ and $W^{\perp}$ mean the tangential part and the normal part of $W\in {\frak X}(N)$, respectively.  
We have, by comparing the tangential part and the normal part of the equation (\ref{eq:3.2}), 
it is equivalent to that 
\begin{align*}
\bigg(\sum_{k=1}^mR^h(\tau(\varphi),d\varphi(e_k))d\varphi(e_k)\bigg)^{\top} &= 0,\quad\mbox{and}\\
\bigg(\sum_{k=1}^mR^h(\tau(\varphi),d\varphi(e_k))\,d\varphi(e_k)\bigg)^{\perp}
&= \sum_{j,k=1}^mh(\tau(\varphi),B_{\varphi}(e_j,e_k))\,B_{\varphi}(e_j,e_k).
\end{align*}
These two equations are equivalent to the following single equation:
\begin{equation}\label{eq:3.3}
\sum_{k=1}^mR^h(\tau(\varphi),d\varphi(e_k))\,d\varphi(e_k)
= \sum_{j,k=1}^mh(\tau(\varphi),B_{\varphi}(e_j,e_k))\,B_{\varphi}(e_j,e_k).
\end{equation}
\vskip0.6cm\par
Summarizing the above, we obtain: 
\begin{theorem} \label{thm:3.2}
Let $\varphi:\,(M^m,g)\rightarrow (N^n,h)$ be an isometric immersion. 
Assume that $\overline{\nabla}^{\perp}_X\tau(\varphi)=0$ for all $X\in{\frak X}(M)$.
Then, $\varphi$ is biharmonic if and only if (\ref{eq:3.3}) holds.
\end{theorem}
\vskip0.6cm\par
As a corollary of Theorem~\ref{thm:3.2}, we obtain:  
\begin{corollary} \label{cor:3.3}
Assume that the sectional curvature of the target space $(N^n,h)$ is non-positive. 
Let 
$\varphi:\,(M^m,g)\rightarrow (N^n,h)$ be an isometric immersion whose tension field satisfies 
$\overline{\nabla}^{\perp}_X\tau(\varphi)=0$ for all $X\in {\frak X}(M)$. Then, if $\varphi$ is biharmonic, then it is harmonic. 
\end{corollary}
\par
{\it Proof.}\quad  
Due to Theorem~\ref{thm:3.2},
taking the inner product to the both hand side of (\ref{eq:3.3}) with $\tau(\varphi)$ which is normal, 
\begin{align}
\sum_{j,k=1}^mh(\tau(\varphi),B_{\varphi}(e_j,e_k))^2
&=h\big(\tau(\varphi),\sum_{k=1}^mR^h(\tau(\varphi),d\varphi(e_k))d\varphi(e_k)\big)\nonumber\\
&=\sum_{k=1}^m h\big(\tau(\varphi),R^h(\tau(\varphi),d\varphi(e_k))d\varphi(e_k)\big)\nonumber\\
&\leq 0. \label{eq:3.4}
\end{align}
Because the quantity 
$h\big(\tau(\varphi),R^h(\tau(\varphi),d\varphi(e_k))d\varphi(e_k)\big)$ is a multiple of the non-negative number
$h(\tau(\varphi),\tau(\varphi))$ times 
the sectional curvature of $(N,h)$ along 
the plane 
$\{\tau(\varphi),d\varphi(e_k)\}$ $(k=1,\cdots,m)$ 
which are also non-positive by our assumption. 
Therefore, the both hand sides of (\ref{eq:3.4}) must be zero, i.e.,
\begin{equation*}
\begin{cases}
h\big(\tau(\varphi), R^h(\tau(\varphi), d\varphi(e_k))\,d\varphi(e_k)\big)=0 & (\forall\, k=1,\cdots,m), \quad \text{and}\\
h(\tau(\varphi), B_{\varphi}(e_j,e_k))=0 & (\forall\, j,k=1,\cdots,m).
\end{cases}
\end{equation*}
Therefore, we obtain 
\begin{equation*}
h(\tau(\varphi),\tau(\varphi))
=h\big(\tau(\varphi),\sum_{j=1}^mB_{\varphi}(e_j,e_j)\big)=0
\end{equation*}
which implies that 
$\tau(\varphi)=0$. 
\qed 
\vskip0.6cm\par
\begin{remark}
Corollary~\ref{cor:3.3} give a partial evidence to the generalized B.-Y. Chen's conjecture:
every biharmonic isometric immersion into a non-positive curvature manifold must be harmonic.
On the other hand, notice that the generalized B.-Y. Chen's conjecture was given by a counter example due to Y. Ou and L. Tang, in 2012, \cite{OT}.
\end{remark}
\vskip0.6cm\par
{\bf 3.2.} In the second part of this section, 
we apply Theorem~\ref{thm:3.2}, to an isometric immersion into an Einstein manifold $(N^n,h)$
whose Ricci transform $\rho^h(u):=\sum_{i=1}^nR^h(u,e'_i)e'_i \ (u\in T_yN,\, y\in N)$, 
where $\{e'_i\}_{i=1}^n$ is a locally defined orthonormal frame field on $(N^n,h)$.    
\par
Let $\{(\xi_i)_x\}_{i=1}^p$ $(x\in M)$ be an orthonormal basis of the orthogonal complement 
$T^{\perp}_{\varphi(x)}M$ of the tangent space $d\varphi_x(T_xM)$ of $M$ in the one 
$T_{\varphi(x)}N$ of $N$ with respect to $h_{\varphi(x)}$: 
$$
T_{\varphi(x)}N=d\varphi_x(T_xM)\oplus \sum_{i=1}^p \mathbb{R}\,(\xi_i)_{x} \qquad (x\in M), 
$$
where 
$p:=\dim N-\dim M=n-m$.
Then, we have
\begin{equation*}
\rho^h(u)=\sum_{k=1}^mR^h(u,d\varphi_k(e_k))d\varphi(e_k)+\sum_{i=1}^pR^h(u,\xi_i)\xi_i
\end{equation*}
for every $u\in T_yN\,(y\in N)$ 
since $\varphi$ is 
an isometric immersion of $(M^m,g)$ into $(N^n,h)$. 
Therefore, in Theorem~\ref{thm:3.2}, we have that 
(\ref{eq:3.3}) is equivalent to the following: 
\begin{equation} \label{eq:3.5}
\rho^h(\tau(\varphi))
-\sum_{i=1}^pR^h(\tau(\varphi),\xi_i)\xi_i
=\sum_{j,k=1}^mh
\big(\tau(\varphi), B_{\varphi}(e_j,e_k)\big)\,
B_{\varphi}(e_j,e_k).
\end{equation}
\par 
Now assume that $(N,h)$ is Einstein, namely, the Ricci transform $\rho^h$ satisfies that 
$\rho^h=c\,\, {\rm I\!d}$ for some constant $c$, 
where ${\rm I\!d}$ is the identity transform. Then, since 
$\rho^h(\tau(\varphi))=c\,\tau(\varphi)$,
we have that (\ref{eq:3.5}) is equivalent to
\begin{equation} \label{eq:3.6}
c\,\tau(\varphi)-\sum_{i=1}^pR^h(\tau(\varphi),\xi_i)\xi_i =\sum_{j,k=1}^mh\big(
\tau(\varphi),B_{\varphi}(e_j,e_k))
\big)\,B_{\varphi}(e_j,e_k).
\end{equation}
Thus we obtain 
\begin{theorem} 
Assume that $\varphi:\,(M^m,g)\rightarrow (N^n,h)$ is an isometric immersion whose tension field 
$\tau(\varphi)$ satisfies that $\overline{\nabla}^{\perp}_X\tau(\varphi)=0$.
\begin{enumerate}
\item Then, $\varphi$ is biharmonic if and only if (\ref{eq:3.5}) holds.
\item In particular, if the target space $(N,h)$ is Einstein, i.e., the 
Ricci transform 
$\rho^h$ of $(N,h)$ satisfies $\rho^h=c\,{\rm I\!d}$ for some constant $c$,
then $\varphi$ is biharmonic if and only if 
(\ref{eq:3.6}) holds. 
\end{enumerate}

\end{theorem}
\medskip\par
{\bf 3.3.}\quad In the following, we treat with a hypersurface
$\varphi:\,(M^m,g)\rightarrow (N^n,h)$, 
i.e., 
$p=1$, and 
$m=\dim M=\dim N-1=n-1$. In this case, let us 
$\xi=\xi_1$ be a unit normal vector field along $\varphi$, and 
denote the second fundamental form $B_{\varphi}$ as 
$B_{\varphi}(e_j,e_k)=H_{jk}\,\xi$ $(j,k=1,\cdots,m)$.
Then, 
\begin{align*}
\lefteqn{\tau(\varphi)=\sum_{j=1}^mB_{\varphi}(e_j,e_j)=\left(\sum_{j=1}^mH_{jj}\right)\,\xi,} \hspace{10mm}\\
\lefteqn{R^h(\tau(\varphi),\xi)\xi=R^h\bigg(\sum_{j=1}^mH_{jj}\,\xi,\xi\bigg)\xi=0,} \hspace{10mm}\\
\lefteqn{\sum_{j,k=1}^mh\big(\tau(\varphi), B_{\varphi}(e_j,e_k))\big)\, B_{\varphi}(e_j,e_k)} \hspace{10mm}\\
&=\big(\sum_{i=1}^mH_{ii}\big)\,\left(\sum_{j,k=1}^mH_{jk}{}^2\right)\,\xi
=\left(\sum_{i=1}^mH_{ii}\right)\,\Vert B_{\varphi}\Vert^2\,\xi,
\end{align*}
where 
$\Vert B_{\varphi}\Vert^2=\sum_{j,k=1}^m\Vert B_{\varphi}(e_j,e_k)\Vert^2=\sum_{j,k=1}^mH_{jk}{}^2$. 
Therefore, 
(\ref{eq:3.5}) holds if and only if 
\begin{equation*}
\left(\sum_{j=1}^mH_{jj}\right) \rho^h(\xi)=\left(\sum_{j=1}^mH_{jj}\right)\,\Vert B_{\varphi}\Vert^2\,\xi
\end{equation*}
which is equivalent to that, 
either 
$\varphi$ is harmonic, i.e., $\sum_{j=1}^mH_{jj}=0$, or 
\begin{equation*}
\rho^h(\xi)=\Vert B_{\varphi}\Vert^2\,\xi.
\end{equation*}
\medskip\par
Thus, we obtain the following theorem: 
\begin{theorem} \label{Theorem 3.6}
Assume that $\varphi:\,(M^m,g)\rightarrow (N^n,h)$ is an isometric immersion whose tension field $\overline{\nabla}_X^{\perp}\tau(\varphi)=0$ $(\forall\,\,X\in {\frak X}(M))$ and $\varphi$ is hypersurface, i.e., $m=n-1$.  
\begin{enumerate}
\item If $\varphi$ is not harmonic, then 
$\varphi$ is biharmonic if and only if 
\begin{equation*}
\rho^h(\xi)=\Vert B_{\varphi}\Vert^2\,\xi,
\end{equation*} 
where $\rho^h$ is the Ricci transform of $(N,h)$, and 
$\xi$ is a unit normal vector field along $\varphi$. 
\item In particular, if $(N,h)$ is an Einstein manifold, i.e., 
$\rho^h=c\,{\rm I\!d}$, and $\varphi$ is not harmonic, then 
$\varphi$ is biharmonic if and only if $\Vert B_{\varphi}\Vert^2=c$.  
\end{enumerate}
\end{theorem}
\medskip\par
Furthermore, we have 
\begin{theorem} \label{thm:3.7}
Assume that
$\varphi:\,(M,g)\rightarrow (N,h)$ is an isometric immersion into a Riemannian manifold $(N,h)$ whose Ricci curvature 
is non-positive,  $\dim M=\dim N-1$, 
and 
$\overline{\nabla}^{\perp}_X\tau(\varphi)=0$ for all $C^{\infty}$ vector field $X$ on $M$. 
Then, if $\varphi$ is biharmonic, it is harmonic. 
\end{theorem}
{\it Proof.}\quad 
If we assume $\varphi$ is not harmonic, then due to (1) of Theorem~\ref{Theorem 3.6}, 
we have 
\begin{equation*}
\rho^h(\xi)=\Vert B_{\varphi}\Vert^2\,\xi.
\end{equation*}
Together with the assumption of non-positivity of the Ricci curvature of $(N,h)$, we have 
\begin{equation*}
0\leq \Vert B_{\varphi}\Vert^2\,h(\xi,\xi)=h(\rho^h(\xi),\xi)\leq 0.
\end{equation*}
Therefore, we have 
\begin{equation*}
h(\rho^h(\xi),\xi)=\Vert B_{\varphi}\Vert^2=0.
\end{equation*}
Hence we have $B_{\varphi}\equiv 0$, 
in particular, we have that 
$\tau(\varphi)=\sum_{i=1}^mB_{\varphi}(e_i,e_i)=0$ which contradicts the assumption. 
\qed 
\medskip\par
Finally, in this section, on the condition 
$\overline{\nabla}^{\perp}_X\tau(\varphi)=0$ $(\forall\,\,X\in {\frak X}(M))$, we give the following criterion:   
\begin{proposition} \label{prop:3.8}
Assume that $\varphi:\,(M^m,g)\rightarrow (N^n,h)$ is an isometric immersion with $m=\dim M=\dim N-1=n-1$. 
Then, the following equivalence holds: The condition that 
$\overline{\nabla}^{\perp}_X\tau(\varphi)=0 \ (\forall\, X\in {\frak X}(M))$ holds
if and only if the mean curvature ${\bf H}=(1/m)\sum_{i=1}^mH_{ii}$ is constant on $M$.
Here, 
$B_{\varphi}(e_i,e_j)=H_{ij}\,\xi$, 
and $\xi$ is a unit normal vector field along $\varphi$. 
\end{proposition}
{\it Proof.}\quad 
Since 
$\tau(\varphi)=m\,{\bf H}\,\xi$ 
where ${\bf H}=(1/m)\sum_{i=1}^mH_{ii}$, 
we have, for all $C^{\infty}$ vector field $X$ on $M$,  
\begin{equation} \label{eq:3.7}
0=\overline{\nabla}^{\perp}_X\tau(\varphi)=
m\,(X{\bf H})\,\xi+m\,{\bf H}\,\overline{\nabla}^{\perp}_X\xi. 
\end{equation}
Since $h(\xi,\xi)=1$, we have 
\begin{align} \label{eq:3.8}
0=\frac12 X_x\,h(\xi,\xi)_{\varphi(x)}
=h({\nabla^h}_{d\varphi(X)}\xi,\xi)
=h(\overline{\nabla}^{\perp}_X\xi,\xi). 
\end{align}
By taking the inner product (\ref{eq:3.7}) and $\xi$, we have, due to (\ref{eq:3.8}),
\begin{equation*}
0=m\,(X\,{\bf H})+m\,{\bf H}\,h(\overline{\nabla}^{\perp}_X\xi,\xi)=
m\,(X\,{\bf H}). 
\end{equation*}
We have $X\,{\bf H}=0$ for all $C^{\infty}$ vector field $X$ on $M$, which implies that ${\bf H}$ is constant on $M$. 
\par 
Conversely, if ${\bf H}$ is constant on $M$, then we have 
\begin{equation*}
\overline{\nabla}^{\perp}_X\tau(\varphi)=m\,(X\,{\bf H})\,\xi+m\,{\bf H}\,\overline{\nabla}^{\perp}_X\xi=m\,{\bf H}\,\overline{\nabla}^{\perp}_X\xi,
\end{equation*}
so that we have, due to (\ref{eq:3.8}), 
\begin{equation*}
h(\overline{\nabla}^{\perp}_X\tau(\varphi),\xi)
=m\,{\bf H}\, h(\overline{\nabla}^{\perp}_X\xi,\xi)=0, 
\end{equation*}
which implies that 
$\overline{\nabla}^{\perp}_X\tau(\varphi)=0$ since it is normal. 
\qed
\vskip0.6cm\par
Summarizing Theorems~\ref{Theorem 3.6} and \ref{thm:3.7}, and Proposition~\ref{prop:3.8}, we obtain
\begin{corollary} \,\,
Let $\varphi:\,(M^m,g)\rightarrow (N^n, h)$ be an isometric immersion. 
Assume that $m=\dim M=\dim N-1=n-1$, and 
the mean curvature of $\varphi$,
${\bf H}=(1/m) \sum_{i=1}^mH_{ii}=(1/m)h(\tau(\varphi),\xi)$, is constant. Then, the following hold: 
\begin{enumerate}
\item Assume that ${\bf H}\not=0$, i.e., $\varphi$ is not harmonic. Then,  it holds that 
$\varphi$ is biharmonic if and only if $\rho^h(\xi)=\Vert B_{\varphi}\Vert^2\,\xi$, where $\rho^h$ is the Ricci transform of $(N,h)$, $\xi$ is a unit normal vector field along $\varphi$ and $B_{\varphi}$ is the second fundamental form of $\varphi$. 
\item Assume that ${\bf H}\not=0$ and 
$(N,h)$ is Einstein, i.e., $\rho^h=c\,\mbox{\rm I\!d}$ for some constant $c$. 
Then, $\varphi$ is biharmonic if and only if $\Vert B_{\varphi}\Vert^2=c$. 
\item Assume that ${\bf H}\not=0$ and the Ricci curvature of $(N,h)$ is non-positive. 
Then, if $\varphi$ is biharmonic, it is harmonic. 
\end{enumerate}
\end{corollary}

\vskip0.6cm\par
\section{Hermann actions and symmetric triads} \label{sect:Hermann_actions}

From this section, we apply the results in Section 3 to the orbits of Hermann actions
using symmetric triads (cf. \cite{Ik1}, \cite{Ik2}, \cite{IST}),
and determine biharmonic regular orbits of cohomogeneity one Hermann actions.
For this purpose,
we express the tension field and the square norm of the second fundamental form
of orbits of Hermann actions in terms of symmetric triads
(Theorem~\ref{thm:norm_of_2nd_fundamental_form1}).

\vskip0.6cm\par
{\bf 4.1.}
First we recall the notions of root system and symmetric triad.
See \cite{Ik1} for details.

Let $(\mathfrak{a}, \langle \cdot, \cdot \rangle)$ be a finite dimensional inner product space over $\mathbb{R}$.
For each $\alpha \in \mathfrak{a}$, we define an orthogonal transformation $s_{\alpha}:\mathfrak{a}\to\mathfrak{a}$ by 
$$
s_{\alpha}(H) = H - \frac{2 \langle \alpha, H \rangle}{\langle \alpha, \alpha \rangle}\alpha \quad (H \in \mathfrak{a}),
$$
namely $s_{\alpha}$ is the reflection with respect to the hyperplane $\{H \in \mathfrak{a} \mid \langle \alpha, H \rangle=0\}$.
\begin{definition}
A finite subset $\Sigma$ of $\mathfrak{a}\setminus \{0\}$ is a root system of $\mathfrak{a}$,
if it satisfies the following three conditions:
\begin{enumerate}
\item $\mathrm{Span}(\Sigma) = \mathfrak{a}$.
\item If $\alpha, \beta \in \Sigma$, then $s_{\alpha}(\beta) \in \Sigma$.
\item $2\langle \alpha, \beta \rangle / \langle \alpha, \alpha \rangle \in \mathbb{Z} \quad (\alpha, \beta \in \Sigma)$.
\end{enumerate}
A root system of $\mathfrak{a}$ is said to be irreducible if it cannot be decomposed into two disjoint nonempty orthogonal subsets.
\end{definition}
Let $\Sigma$ be a root system of $\mathfrak{a}$. 
The Weyl group $W(\Sigma)$ of $\Sigma$ is the finite subgroup of the orthogonal group $\mathrm{O}(\mathfrak{a})$ of $\mathfrak{a}$
generated by $\{ s_{\alpha} \mid \alpha \in \Sigma\}$.
\begin{definition}[\cite{Ik1} Definition 2.2]
A triple $(\tilde{\Sigma}, \Sigma, W)$ of finite subsets of $\mathfrak{a} \setminus \{0\}$ is a symmetric triad of $\mathfrak{a}$, 
if it satisfies the following six conditions:
\begin{enumerate}
\item $\tilde{\Sigma}$ is an irreducible root system of $\mathfrak{a}$. 
\item $\Sigma$ is a root system of $\mathfrak{a}$.
\item $(-1)W=W, \ \tilde{\Sigma}=\Sigma\cup W$.
\item $\Sigma \cap W$ is a nonempty subset.
If we put $l:=\max\{\|\alpha \| \mid \alpha \in \Sigma \cap W \}$,
then $\Sigma \cap W = \{\alpha \in \tilde{\Sigma} \mid \|\alpha \| \leq l\}$.
\item For $\alpha \in W$ and $\lambda \in \Sigma \setminus W$,
$$
2\frac{\langle \alpha, \lambda \rangle}{\langle \alpha, \alpha \rangle} \ \text{is odd if and only if} \ s_{\alpha}(\lambda) \in W \setminus \Sigma.
$$
\item For $\alpha \in W$ and $\lambda \in W \setminus \Sigma$,
$$
2\frac{\langle \alpha, \lambda \rangle}{\langle \alpha, \alpha \rangle} \ \text{is odd if and only if } \ s_{\alpha}(\lambda) \in \Sigma \setminus W.
$$ 
\end{enumerate}
\end{definition}
We define an open subset $\mathfrak{a}_r$ of $\mathfrak{a}$ by
$$
\mathfrak{a}_{r} = \bigcap_{\lambda \in \Sigma, \alpha \in W} \left\{H \in \mathfrak{a} \ \Big| \
\langle \lambda, H \rangle \not\in \pi \mathbb{Z},\ \langle \alpha, H \rangle \not\in \frac{\pi}{2} + \pi \mathbb{Z} \right\}.
$$
A point in $\mathfrak{a}_{r}$ is called a regular point,
and a point in the complement of $\mathfrak{a}_{r}$ in $\mathfrak{a}$ is called a singular point.
A connected component of $\mathfrak{a}_{r}$ is called a cell.
The affine Weyl group $\tilde{W}(\tilde{\Sigma}, \Sigma, W)$ of $(\tilde{\Sigma}, \Sigma, W)$ is a subgroup of the affine group of $\mathfrak{a}$,
i.e. the semidirect product $\mathrm{O}(\mathfrak{a}) \ltimes \mathfrak{a}$, generated by 
$$
\left\{\left(s_{\lambda}, \frac{2n\pi}{\| \lambda\|^2}\lambda\right) \ \bigg| \ \lambda \in \Sigma, n \in \mathbb{Z}\right\}
\cup \left\{\left(s_{\alpha}, \frac{(2n+1)\pi}{\| \alpha\|^2}\alpha\right) \ \bigg| \ \alpha \in W, n \in \mathbb{Z}\right\}.
$$
The action of $\left(s_{\lambda}, (2n\pi/\| \lambda\|^2)\lambda\right)$ on $\mathfrak{a}$
is the reflection with respect to the hyperplane $\{ H \in \mathfrak{a} \mid \langle \lambda, H \rangle = n\pi \}$, 
and the action of $\left(s_{\alpha}, ((2n+1)\pi/\| \alpha\|^2)\alpha\right)$ on $\mathfrak{a}$
is the reflection with respect to the hyperplane $\{ H \in \mathfrak{a} \mid \langle \alpha, H \rangle = (n+1/2)\pi\}$.
The affine Weyl group $\tilde{W}(\tilde{\Sigma}, \Sigma, W)$ acts transitively on the set of all cells.
More precisely, for each cell $P$, it holds that
$$
\mathfrak{a} = \bigcup_{s \in \tilde{W}(\tilde{\Sigma}, \Sigma, W)} s \overline{P}.
$$

We take a fundamental system $\tilde{\Pi}$ of $\tilde{\Sigma}$.
We denote by $\tilde{\Sigma}^{+}$ the set of positive roots in $\tilde{\Sigma}$.
Set $\Sigma^{+} = \tilde{\Sigma}^{+} \cap \Sigma$ and $W^{+} = \tilde{\Sigma}^{+} \cap W$.
Denote by $\Pi$ the set of simple roots of $\Sigma$.
We set 
$$
W_{0} = \{\alpha \in W^{+} \mid \alpha + \lambda \not\in W \ (\lambda \in \Pi) \}.
$$
From the classification of symmetric triads,
we have that $W_{0}$ consists of only one element, denoted by $\tilde{\alpha}$. 
We define an open subset $P_0$ of $\mathfrak{a}$ by
\begin{equation} \label{eq:cell}
P_{0} = \left\{H \in \mathfrak{a} \ \Big| \ \langle \tilde{\alpha}, H \rangle < \frac{\pi}{2},\
\langle \lambda, H \rangle >0 \ (\lambda \in \Pi) \right\}.
\end{equation}
Then $P_{0}$ is a cell. 
\begin{definition}[\cite{Ik1} Definition 2.13]
Let $(\tilde{\Sigma}, \Sigma, W)$ be a symmetric triad of $\mathfrak{a}$.
Consider two mappings $m$ and $n$ from $\tilde{\Sigma}$ to $\mathbb{R}_{\geq 0} := \{a \in \mathbb{R} \mid a \geq 0\}$
which satisfy the following four conditions:
\begin{enumerate}
\item For any $\lambda \in \tilde{\Sigma}$,
  \begin{enumerate}
  \item[(1-1)] $m(\lambda) = m(-\lambda)$,\ $n(\lambda) = n(-\lambda)$,
  \item[(1-2)] $m(\lambda) > 0$ if and only if $\lambda \in \Sigma$,
  \item[(1-3)] $n(\lambda) > 0$ if and only if $\lambda \in W$.
  \end{enumerate}
\item When $\lambda \in \Sigma, \alpha \in W, s \in W(\Sigma)$,
then $m(\lambda) = m(s(\lambda)),\ n(\alpha)=n(s(\alpha))$.
\item When $\lambda \in \tilde{\Sigma},\ \sigma \in W(\tilde{\Sigma})$,
then $m(\lambda) + n(\lambda) = m(\sigma(\lambda)) + n(\sigma(\lambda))$.
\item Let $\lambda \in \Sigma \cap W$, $\alpha \in W$. 
If $2\langle \alpha, \lambda \rangle / \langle \alpha, \alpha \rangle$ is even,
then $m(\lambda)= m(s_{\alpha}(\lambda))$.
If $2\langle \alpha, \lambda \rangle / \langle \alpha, \alpha \rangle$ is odd,
then $m(\lambda)= n(s_{\alpha}(\lambda))$.
\end{enumerate}
We call $m(\lambda)$ and $n(\alpha)$ the multiplicities of $\lambda$ and $\alpha$, respectively.
\end{definition}

\vskip0.6cm\par
{\bf 4.2.}
We will review some basics of the theory of compact symmetric spaces.

Let $G$ be a compact connected Lie group and $K$ a closed subgroup of $G$.
Assume that there exists an involutive automorphism $\theta$ of $G$ which satisfies
$(G_{\theta})_0 \subset K \subset G_{\theta}$,
where $G_{\theta}$ is the set of fixed points of $\theta$ and $(G_{\theta})_0$ is the identity component of $G_{\theta}$.
Then the pair $(G, K)$ is called a compact symmetric pair.
We denote the Lie algebras of $G$ and $K$ by $\mathfrak{g}$ and $\mathfrak{k}$, respectively.
The involutive automorphism $\theta$ of $G$ induces an involutive automorphism of $\mathfrak{g}$,
which is also denoted by the same symbol $\theta$.
We can see that
$$
\mathfrak{k} = \{X \in \mathfrak{g} \mid \theta(X) = X \},
$$
and we define
$$
\mathfrak{m} = \{X \in \mathfrak{g} \mid \theta(X) = -X \}.
$$
Take an inner product $\langle \cdot, \cdot \rangle$ on $\mathfrak{g}$
which is invariant under the actions of $\mathrm{Ad}(G)$ and $\theta$.
The inner product $\langle \cdot, \cdot \rangle$ induces a bi-invariant Riemannian metric on $G$
and a $G$-invariant Riemannian metric on $N=G/K$, which are denoted by the same symbol $\langle \cdot, \cdot \rangle$.
Then $(N, \langle \cdot, \cdot \rangle)$ is a compact symmetric space.
Conversely, any compact symmetric space can be constructed in this way.
Since $\theta$ is involutive,
we have an orthogonal direct sum decomposition of $\mathfrak{g}$:
$$
\mathfrak{g} = \mathfrak{k} \oplus \mathfrak{m}.
$$
This decomposition is called the canonical decomposition of $(G, K)$.
We denote by $\pi$ the natural projection from $G$ onto $N$. 
The tangent space $T_{\pi(e)}N$ of $N$ at the origin $\pi(e)$ is identified with $\mathfrak{m}$ in a natural way.
The Ricci tensor $\mathrm{Ric}(\cdot, \cdot)$ of $N$ is given by
$$
\mathrm{Ric}(X,Y) = -\frac{1}{2} \mathrm{Killing}(X,Y) \quad (X, Y \in \mathfrak{m}),
$$
where $\mathrm{Killing}(\cdot ,\cdot)$ is the Killing form of $\mathfrak{g}$.
If $G$ is semisimple, then we can give an $\mathrm{Ad}(G)$-invariant inner product on $\mathfrak{g}$
by $\langle \cdot, \cdot \rangle = -\mathrm{Killing}(\cdot, \cdot)$,
hence $N$ is an Einstein manifold with Einstein constant $c = 1/2$. 

Here, let us recall the notion of hyperpolar actions (cf. \cite{Kollross}).
An isometric action of a compact Lie group on a Riemannian manifold is said to be hyperpolar
if there exists a closed, connected submanifold that is flat in the induced metric and meets all orbits orthogonally.
Such a submanifold is called a section of the Lie group action.
Kollross \cite{Kollross} classified hyperpolar actions on irreducible symmetric spaces of compact type.

\vskip0.6cm\par
{\bf 4.3.}
Our aim is to apply the theory of symmetric triad due to Ikawa \cite{Ik1}
in order to express the second fundamental form of orbits of Hermann actions.

Let $(G, K_1)$ and $(G, K_2)$ be compact symmetric pairs with respect to involutive automorphisms $\theta_1$ and $\theta_2$ of a compact Lie group $G$, respectively.
Then the triple $(G,K_{1},K_{2})$ is called a compact symmetric triad.
We denote the Lie algebras of $G, K_{1}$ and $K_{2}$ by $\mathfrak{g}, \mathfrak{k}_{1}$ and $\mathfrak{k}_{2}$, respectively.
The involutive automorphism of $\mathfrak{g}$ induced from $\theta_{i}$ will be also denoted by $\theta_{i}$.
Take an $\mathrm{Ad}(G)$-invariant inner product $\langle \cdot, \cdot \rangle$ on $\mathfrak{g}$.
Then the inner product $\langle \cdot, \cdot \rangle$ induces a bi-invariant Riemannian metric on $G$
and $G$-invariant Riemannian metrics on the coset manifolds $N_{1}=G/K_{1}$ and $N_2=G/K_2$.
We denote these Riemannian metrics on $G$, $N_1$ and $N_2$ by the same symbol $\langle \cdot, \cdot \rangle$.
The isometric action of $K_{2}$ on $N_{1}$ and the action of $K_1$ on $N_2$ are called Hermann actions.
Now we have two canonical decompositions of $\mathfrak{g}$:
$$
\mathfrak{g} = \mathfrak{k}_{1} \oplus \mathfrak{m}_{1} = \mathfrak{k}_{2} \oplus \mathfrak{m}_{2},
$$
where $\mathfrak{m}_{i} = \{ X \in \mathfrak{g} \mid \theta_{i}(X) =-X \}\ (i=1,2)$.
We define a closed subgroup $G_{12}$ of $G$ by
$$
G_{12}=\{ k\in G \mid \theta_{1} (k) = \theta_{2} (k)\},
$$
and we denote the identity component of $G_{12}$ by $(G_{12})_{0}$.
Then $\theta_1$ induces an involutive automorphism of $(G_{12})_{0}$.
Hence $((G_{12})_0, K_{12})$ is a compact symmetric pair,
where $K_{12}$ is a closed subgroup of $(G_{12})_{0}$ defined by 
$$
K_{12} = \{k \in (G_{12})_{0} \mid \theta_{1}(k) = k\}.
$$
The canonical decomposition of the Lie algebra $\mathfrak{g}_{12}$ of $(G_{12})_{0}$ is given by 
$$
\mathfrak{g}_{12} = (\mathfrak{k}_{1} \cap \mathfrak{k}_{2}) \oplus (\mathfrak{m}_{1} \cap \mathfrak{m}_{2}).
$$
Fix a maximal abelian subspace $\mathfrak{a}$ in $\mathfrak{m}_{1} \cap \mathfrak{m}_{2}$. 
Then $\exp \mathfrak{a}$ is a torus subgroup in $(G_{12})_{0}$.
We denote by $\pi_{i}$ the natural projection from $G$ onto $N_{i} \, (i=1,2)$.
Then, the totally geodesic flat torus $\pi_{1}(\exp \mathfrak{a})$ is a section of $K_{2}$-action on $N_{1}$,
hence the action is hyperpolar.
Similarly $\pi_{2}(\exp \mathfrak{a})$ is a section of $K_{1}$-action on $N_{2}$.
The cohomogeneity of $K_{2}$-action on $N_{1}$ and that of $K_1$-action on $N_2$ are equal to $\dim \mathfrak{a}$.
We call an orbit of the maximal dimension a regular orbit.
For $k \in G$, we denote the left transformation of $G$ by $L_{k}$.
The isometries on $N_{1}$ and $N_2$ induced by $L_{k}$ will be also denoted by the same symbol $L_{k}$.

We should prepare several terminologies to determine the second fundamental forms
of the regular orbits of Hermann actions.
\begin{definition}
A compact symmetric triad $(G, K_1, K_2)$ is said to be commutative
if $\theta_1\theta_2 = \theta_2\theta_1$.
Then $K_2$-action on $N_1$ and $K_1$-action on $N_2$ are called
commutative Hermann actions.
\end{definition}
Hereafter we assume that $(G, K_1, K_2)$ is a commutative compact symmetric triad
where $G$ is semisimple.
Then we have 
$$
\mathfrak{g} = (\mathfrak{k}_{1} \cap \mathfrak{k}_{2}) \oplus (\mathfrak{m}_{1} \cap \mathfrak{m}_{2})
\oplus (\mathfrak{k}_{1} \cap \mathfrak{m}_{2}) \oplus (\mathfrak{m}_{1} \cap \mathfrak{k}_{2}).
$$
We define subspaces in $\mathfrak{g}$ as follows:
\begin{align*}
\mathfrak{k}_{0} &= \{ X \in \mathfrak{k}_{1}\cap \mathfrak{k}_{2} \mid [\mathfrak{a}, X] =\{0\}\},\\
V(\mathfrak{k}_{1} \cap \mathfrak{m}_{2}) &= \{ X \in \mathfrak{k}_{1} \cap \mathfrak{m}_{2} \mid [\mathfrak{a}, X] =\{0\}\},\\  
V(\mathfrak{m}_{1} \cap \mathfrak{k}_{2}) &= \{ X \in \mathfrak{m}_{1} \cap \mathfrak{k}_{2} \mid [\mathfrak{a}, X] =\{0\}\},
\end{align*}
and for $\lambda \in \mathfrak{a}$
\begin{align*}
\mathfrak{k}_{\lambda} &= \{X \in \mathfrak{k}_{1}\cap \mathfrak{k}_{2} \mid [H,[H,X]] =-\langle \lambda , H \rangle^{2} X \ (H \in \mathfrak{a})\},\\
\mathfrak{m}_{\lambda} &= \{X \in \mathfrak{m}_{1}\cap \mathfrak{m}_{2} \mid [H,[H,X]] =-\langle \lambda , H \rangle^{2} X \ (H \in \mathfrak{a})\},\\
V^{\perp}_{\lambda}(\mathfrak{k}_{1} \cap \mathfrak{m}_{2})
&= \{ X \in \mathfrak{k}_{1} \cap \mathfrak{m}_{2} \mid [H,[H,X]] =-\langle \lambda, H \rangle^{2} X \ (H \in \mathfrak{a})\},\\
V^{\perp}_{\lambda}(\mathfrak{m}_{1} \cap \mathfrak{k}_{2})
&= \{ X \in \mathfrak{m}_{1} \cap \mathfrak{k}_{2} \mid [H,[H,X]] =-\langle \lambda, H \rangle^{2} X \ (H \in \mathfrak{a})\}.
\end{align*}
We set 
\begin{align*}
\Sigma &= \{\lambda \in \mathfrak{a} \setminus \{0\} \mid \mathfrak{k}_{\lambda} \neq \{0\}\},\\ 
W &= \{\alpha \in \mathfrak{a} \setminus \{0\} \mid V^{\perp}_{\alpha}(\mathfrak{k}_{1} \cap \mathfrak{m}_{2}) \neq \{0\}\},\\ 
\tilde{\Sigma} &= \Sigma \cup W.
\end{align*}
It is known that $\dim \mathfrak{k}_{\lambda} = \dim \mathfrak{m}_{\lambda}$
and $\dim V^{\perp}_{\lambda}(\mathfrak{k}_{1} \cap \mathfrak{m}_{2}) = \dim V^{\perp}_{\lambda}(\mathfrak{m}_{1} \cap \mathfrak{k}_{2})$
for each $\lambda \in \tilde{\Sigma}$.
Thus we define 
$m(\lambda) := \dim \mathfrak{k}_{\lambda}$ and $n(\lambda) := \dim V^{\perp}_{\lambda }(\mathfrak{k}_{1} \cap \mathfrak{m}_{2})$.
Notice that $\Sigma$ is the root system of the symmetric pair $((G_{12})_{0}, K_{12})$ with respect to $\mathfrak{a}$. 
\begin{proposition}[\cite{Ik1} Lemma 4.12, Theorem 4.33] \label{prop:compact_symmetric_triad}
Let $(G, K_1, K_2)$ be a commutative compact symmetric triad where $G$ is semisimple.
Then $\tilde{\Sigma}$ is a root system of $\mathfrak{a}$.
In addition, if $G$ is simple and $\theta_1 \not\sim \theta_2$,
then $(\tilde{\Sigma}, \Sigma, W)$ is a symmetric triad of $\mathfrak{a}$,
moreover $m(\lambda)$ and $n(\alpha)$ are multiplicities of $\lambda \in \Sigma$ and $\alpha \in W$.
Here $\theta_1 \not\sim \theta_2$ means that 
$\theta_1$ and $\theta_2$ cannot be transformed each other by an inner automorphism of $\mathfrak{g}$.
\end{proposition}
We take a basis of $\mathfrak{a}$ and define a lexicographic ordering $>$ on $\mathfrak{a}$ with respect to the basis.
We set 
$$
\tilde{\Sigma}^{+} = \{\lambda \in \tilde{\Sigma} \mid \lambda > 0\},\quad
\Sigma^{+} = \Sigma \cap \tilde{\Sigma}^{+},\quad
W^{+} = W \cap \tilde{\Sigma}^{+}.
$$
Then we have an orthogonal direct sum decomposition of $\mathfrak{g}$:
\begin{align*}
\mathfrak{g} = \mathfrak{k}_{0} \oplus \sum_{\lambda \in \Sigma^{+}} \mathfrak{k}_{\lambda}
\oplus \mathfrak{a} \oplus \sum_{\lambda \in \Sigma^{+}} \mathfrak{m}_{\lambda}
&\oplus V(\mathfrak{k}_{1} \cap \mathfrak{m}_{2}) \oplus \sum_{\alpha \in W^{+}} V^{\perp}_{\alpha}(\mathfrak{k}_{1} \cap \mathfrak{m}_{2})\\
&\oplus V(\mathfrak{m}_{1} \cap \mathfrak{k}_{2}) \oplus \sum_{\alpha \in W^{+}} V^{\perp}_{\alpha}(\mathfrak{m}_{1} \cap \mathfrak{k}_{2}).
\end{align*}
Furthermore we have the following lemma.
\begin{lemma}[\cite{Ik1} Lemmas 4.3, 4.16] \label{onb}
\begin{enumerate}
\item For each $\lambda \in \Sigma^{+}$,
there exist orthonormal bases $\{S_{\lambda, i}\}_{i=1}^{m(\lambda)}$
and $\{T_{\lambda, i}\}_{i=1}^{m(\lambda)}$ of $\mathfrak{k}_{\lambda}$ and $\mathfrak{m}_{\lambda}$ respectively such that 
for any $H \in \mathfrak{a}$
$$
[H, S_{\lambda, i}] = \langle \lambda, H \rangle T_{\lambda, i},\quad
[H, T_{\lambda, i}] = -\langle \lambda, H \rangle S_{\lambda, i},\quad
[S_{\lambda, i}, T_{\lambda, i}] = \lambda,
$$
\begin{align*}
\mathrm{Ad}(\exp H) S_{\lambda, i} &= \cos \langle \lambda, H \rangle S_{\lambda, i} + \sin \langle \lambda, H \rangle T_{\lambda, i},\\
\mathrm{Ad}(\exp H) T_{\lambda, i} &= -\sin \langle \lambda, H \rangle S_{\lambda, i} + \cos \langle \lambda, H \rangle T_{\lambda, i}.
\end{align*}
\item For each $\alpha \in W^{+}$,
there exist orthonormal bases 
$\{X_{\alpha, j}\}_{j=1}^{n(\alpha)}$ and $\{Y_{\alpha, j}\}_{j=1}^{n(\alpha)}$
of $V^{\perp}_{\alpha}(\mathfrak{k}_{1} \cap \mathfrak{m}_{2})$
and $V^{\perp}_{\alpha}(\mathfrak{m}_{1} \cap \mathfrak{k}_{2})$ respectively such that for any $H \in \mathfrak{a}$
$$
[H, X_{\alpha, j}] = \langle \alpha, H \rangle Y_{\alpha, j},\quad
[H, Y_{\alpha, j}] = -\langle \alpha, H \rangle X_{\alpha, j},\quad
[X_{\alpha, j}, Y_{\alpha, j}] = \alpha,
$$
\begin{align*}
\mathrm{Ad}(\exp H) X_{\alpha, j} &= \cos \langle \alpha, H \rangle X_{\alpha, j} + \sin \langle \alpha, H \rangle Y_{\alpha, j},\\
\mathrm{Ad}(\exp H) Y_{\alpha, j} &= -\sin \langle \alpha, H \rangle X_{\alpha, j} + \cos \langle \alpha, H \rangle Y_{\alpha, j}.
\end{align*}
\end{enumerate}
\end{lemma}

Now we consider the second fundamental form of an orbit $K_{2} \pi_{1}(x)$ of
the action of $K_{2}$ on $N_{1} = G/K_1$ for $x \in G$ (cf. Lemma~\ref{2nd}).
Without loss of generalities we can assume that $x = \exp H$ where $H \in \mathfrak{a}$,
since $\pi_{1}(\exp \mathfrak{a})$ is a section of the action.
We identify the tangent space $T_{\pi_{1}(e)} N_{1}$ with $\mathfrak{m}_{1}$ via $(d\pi_{1})_{e}$.
For $x= \exp H \ (H \in \mathfrak{a})$,
the tangent space and the normal space of $K_2\pi_1(x)$ at $\pi_1(x)$ are given as
\begin{align*}
dL_{x}^{-1}(T_{\pi_{1}(x)} (K_{2}\pi_{1}(x)))
&\cong (\mathrm{Ad}(x^{-1}) \mathfrak{k}_{2})_{\mathfrak{m}_1}\\
&= \sum_{\lambda \in \Sigma^{+} \atop \langle \lambda, H \rangle \not\in \pi \mathbb{Z}} \mathfrak{m}_{\lambda}
\oplus V(\mathfrak{m}_{1} \cap \mathfrak{k}_{2})
\oplus \sum_{\alpha \in W^{+} \atop \langle \alpha, H \rangle \not\in (\pi/2) + \pi \mathbb{Z}}
V^{\perp}_{\alpha}(\mathfrak{m}_{1} \cap \mathfrak{k}_{2}),\\
dL_{x}^{-1}(T_{\pi_{1}(x)}^{\perp} (K_{2}\pi_{1}(x)))
&\cong (\mathrm{Ad}(x^{-1}) \mathfrak{m}_{2}) \cap \mathfrak{m}_{1}\\
&= \mathfrak{a} \oplus \sum_{\lambda \in \Sigma^{+} \atop \langle \lambda, H \rangle \in \pi \mathbb{Z}} \mathfrak{m}_{\lambda}
\oplus \sum_{\alpha \in W^{+} \atop \langle \alpha, H \rangle \in (\pi/2) + \pi \mathbb{Z}}
V^{\perp}_{\alpha}(\mathfrak{m}_{1} \cap \mathfrak{k}_{2}),
\end{align*}
where $X_{\mathfrak{m}_1}$ denotes $\mathfrak{m}_1$-component of $X \in \mathfrak{g}$ with respect to
the canonical decomposition $\mathfrak{g} = \mathfrak{k}_1 \oplus \mathfrak{m}_1$.
Using the above decompositions of the tangent space and the normal space of the orbit $K_2\pi_1(x)$
and the orthonormal basis given in Lemma~\ref{onb},
we can apply Ikawa's results (cf. Lemma 4.22 in \cite{Ik1}) to our cases.
Let us denote the second fundamental form and the tension field of the orbit $K_{2}\pi_{1}(x)$ in $N_{1}$
by $B_H$ and $\tau_H$, respectively. Then we have the following lemma.
\begin{lemma} \label{2nd}
Let $x=\exp H$ for $H \in \mathfrak{a}$.
Then we have:
\begin{enumerate}
\item $dL_{x}^{-1}B_{H}(dL_{x}(T_{\lambda, i}), dL_{x}(T_{\mu, j})) = \cot(\langle \mu, H \rangle) [T_{\lambda, i}, S_{\mu, j}]^{\perp}$,
\item $dL_{x}^{-1}B_{H}(dL_{x}(Y_{\alpha, i}), dL_{x}(Y_{\beta, j})) = -\tan(\langle \beta, H \rangle) [Y_{\alpha, i}, X_{\beta, j}]^{\perp}$,
\item $B_{H}(dL_{x}(Y_{1}), dL_{x}(Y_{2}))=0$,
\item $B_{H}(dL_{x}(T_{\lambda, i}), dL_{x}(Y_2)) = 0$,
\item $B_{H}(dL_{x}(Y_{\alpha, i}), dL_{x}(Y_2)) = 0$,
\item $dL_{x}^{-1}B_{H}(dL_{x}(T_{\lambda, i}), dL_{x}(Y_{\beta, j})) = -\tan(\langle \beta, H \rangle) [T_{\lambda, i}, X_{\beta, j}]^{\perp}$,
\end{enumerate}
for
\begin{align*}
&\lambda, \mu \in \Sigma^{+} \ \text{with} \ \langle \lambda, H \rangle, \langle \mu, H \rangle \not\in \pi \mathbb{Z},\ 1 \leq i \leq m(\lambda),\ 1 \leq j \leq m(\mu),\\
&\alpha, \beta \in W^{+} \ \text{with} \ \langle \alpha, H \rangle, \langle \beta, H \rangle \not\in \frac{\pi}{2} + \pi \mathbb{Z},\ 1 \leq i \leq n(\alpha),\ 1 \leq j \leq n(\beta),\\
&Y_{1}, Y_{2} \in V(\mathfrak{m}_{1} \cap \mathfrak{k}_{2}).
\end{align*}
Here $X^{\perp}$ is the normal component, i.e. $(\mathrm{Ad}(x^{-1}) \mathfrak{m}_{2}) \cap \mathfrak{m}_{1}$-component,
of a tangent vector $X \in \mathfrak{m}_{1}$. 
\end{lemma}

\par
Due to Lemma~\ref{2nd}, we have the following.
\begin{theorem}\label{thm:norm_of_2nd_fundamental_form1}
If $K_{2}\pi_{1}(x)$ is a regular orbit, then
\begin{align}
&\|B_{H}\|^{2} = \sum_{\lambda \in \Sigma^{+}} m(\lambda) (\cot \langle \lambda, H \rangle)^{2} \langle \lambda, \lambda \rangle
+ \sum_{\alpha \in W^{+}} n(\alpha) (\tan \langle \alpha, H \rangle)^{2} \langle \alpha, \alpha \rangle, \label{biheq} \\
&dL_x^{-1}(\tau_{H}) = -\sum_{\lambda \in \Sigma^{+}} m(\lambda) \cot \langle \lambda, H \rangle \lambda
+ \sum_{\alpha \in W^{+}} n(\alpha) \tan \langle \alpha, H \rangle \alpha.
\label{harmeq}
\end{align}
\end{theorem}
\par
{\it Proof.}\quad
For (\ref{biheq}), since the orbit $K_{2}\pi_{1}(x)$ is regular,
its tangent space and normal space are given as:
\begin{align*}
dL_{x}^{-1} (T_{\pi_{1}(x)} (K_{2}\pi_{1}(x)))
&= \sum_{\lambda \in \Sigma^{+}} \mathfrak{m}_{\lambda} \oplus V(\mathfrak{m}_{1} \cap \mathfrak{k}_{2})
\oplus \sum_{\alpha \in W^{+}} V_{\alpha}^{\perp}(\mathfrak{m}_{1} \cap \mathfrak{k}_{2}), \\
dL_{x}^{-1} (T_{\pi_{1}(x)}^{\perp} (K_{2}\pi_{1}(x))) &= \mathfrak{a}.
\end{align*}
For each $\lambda, \mu \in \Sigma^{+}, 1 \leq i \leq m(\lambda), 1 \leq j \leq m(\mu)$, we have 
\begin{equation*}
\langle [T_{\lambda, i}, S_{\mu, j}], H' \rangle
= \langle T_{\lambda, i}, [S_{\mu, j}, H'] \rangle
= \langle T_{\lambda, i}, -\langle \mu, H' \rangle T_{\mu, j}\rangle 
= -\delta_{\lambda, \mu} \delta_{i, j} \langle \mu, H' \rangle
\end{equation*}
for all $H' \in \mathfrak{a}$.
Thus, we have 
$$
[T_{\lambda ,i}, S_{\mu , j}]^{\perp}= -\delta_{\lambda ,\mu }\delta_{i,j} \mu .
$$
Similarly, we have 
\begin{align*}
[Y_{\alpha, i}, X_{\beta, j}]^{\perp} &= -\delta_{\alpha, \beta} \delta_{i, j} \beta \quad
(\alpha, \beta \in W^{+}, 1 \leq i \leq n(\alpha), 1 \leq j \leq n(\beta)), \\
[T_{\lambda, i}, X_{\beta, j}]^{\perp} &= 0 \quad (\lambda \in \Sigma^{+}, \beta \in W^{+}, 1\leq i \leq m(\lambda), 1 \leq j \leq n(\beta)).
\end{align*}
Form Lemma~\ref{2nd}, we obtain 
\begin{align*}
\|B_{H}\|^2
&= \sum_{\lambda, \mu, i, j} \|\cot \langle \mu, H \rangle [T_{\lambda, i}, S_{\mu, j}]^{\perp} \|^2
+ \sum_{\alpha, \beta, i, j} \|-\tan \langle \beta, H \rangle [Y_{\alpha, i}, X_{\beta, j}]^{\perp} \|^2\\
&= \sum_{\lambda, \mu, i, j} \|-\cot \langle \mu, H \rangle \delta_{\lambda, \mu} \delta_{i, j} \mu \|^2
+ \sum_{\alpha, \beta, i, j} \|\tan \langle \beta, H \rangle \delta_{\alpha, \beta} \delta_{i, j} \beta \|^2\\
&= \sum_{\lambda \in \Sigma^{+}} \sum_{i=1}^{m(\lambda)} (\cot \langle \lambda ,H \rangle)^{2} \langle \lambda, \lambda \rangle 
+\sum_{\alpha  \in W^{+}} \sum_{i=1}^{n(\alpha )} (\tan  \langle \alpha ,H \rangle )^{2}\langle \alpha , \alpha \rangle \\
&= \sum_{\lambda \in \Sigma^{+}} m(\lambda) (\cot \langle \lambda ,H \rangle)^{2} \langle \lambda, \lambda \rangle 
+\sum_{\alpha  \in W^{+}} n(\alpha ) (\tan  \langle \alpha ,H \rangle )^{2}\langle \alpha , \alpha \rangle.
\end{align*}
The formula (\ref{harmeq}) was proved in \cite[Corollary 4.23]{Ik1}.
\qed

\vskip0.6cm\par

\section{Biharmonic orbits of cohomogeneity one Hermann actions} \label{sect:biharmonic_orbits}
In this section, applying Theorem~\ref{Theorem 3.6} we will study biharmonic regular orbits of cohomogeneity one Hermann actions.

Let $(G,K_{1},K_{2})$ be a commutative compact symmetric triad
where $G$ is semisimple.
Define an inner product $\langle \cdot, \cdot \rangle$ on $\mathfrak{g}$ by $\langle \cdot, \cdot \rangle = -\mathrm{Killing}(\cdot, \cdot)$.
Then, $(N_{1}, \langle \cdot, \cdot \rangle)$ and $(N_{2}, \langle \cdot, \cdot \rangle)$ are Einstein manifolds with Einstein constant $c = 1/2$.
It is known that the tension field of an orbit of a Hermann action is parallel in the normal bundle (see \cite{IST}),
i.e. $\overline{\nabla}_X^{\perp}\tau_{H}=0$ for every vector field $X$ on the orbit $K_{2}\pi_{1}(x)$.

Hereafter we assume that $\dim \mathfrak{a} = 1$.
Since the cohomogeneity of $K_{2}$-action on $N_{1}$ and that of $K_1$-action on $N_2$ are equal to $\dim \mathfrak{a}$,
regular orbits of $K_2$-actions (resp. $K_1$-action) are homogeneous hypersurfaces in $N_1$ (resp. $N_2$).
Hence we can apply (2) of Theorem~\ref{Theorem 3.6} for regular orbits of these actions.
Clearly, $K_{2}\pi_{1}(x)$ is a regular orbit if and only if $K_1 \pi_{2}(x)$ is also a regular orbit. 
Therefore, we have the following proposition.
\begin{proposition}\label{prop5.1}
Let $x = \exp H$ for $H \in \mathfrak{a}$.
Suppose that $K_{2}\pi_{1}(x)$ is a regular orbit of $K_{2}$-action on $N_{1}$,
so $K_1\pi_{2}(x)$ is also a regular orbit of $K_{1}$-action on $N_{2}$.
Then,
\begin{enumerate}
\item An orbit $K_{2}\pi_{1}(x)$ is harmonic if and only if $K_1\pi_{2}(x)$ is harmonic.
\item An orbit $K_{2}\pi_{1}(x)$ is proper biharmonic if and only if $K_1\pi_{2}(x)$ is proper biharmonic.
\end{enumerate}
\end{proposition}

\par
{\it Proof.}\quad
Analogous to Lemma~\ref{2nd}, we can express the second fundamental form $B'_H$ of
$K_1 \pi_2(x)$ in $N_2$ using the orthonormal basis given in Lemma~\ref{onb}.
Then easily we can verify
$$
\|B'_H\|^2 = \|B_H\|^2,\qquad
dL_x^{-1} (\tau'_H) = dL_x^{-1} (\tau_H),
$$
where $\tau'_H$ denotes the tension field of $K_1 \pi_2(x)$ in $N_2$.
Therefore, from Theorem~\ref{Theorem 3.6}, we have the consequence.
\qed 

\vskip0.6cm\par
If $G$ is simple and $\theta_1 \not\sim \theta_2$,
then for a commutative compact symmetric triad $(G, K_1, K_2)$ the triple $(\tilde{\Sigma}, \Sigma, W)$ is a symmetric triad
with multiplicities $m(\lambda)$ and $n(\alpha)$ (cf. Proposition~\ref{prop:compact_symmetric_triad}).
In this case, for $x = \exp H \ (H \in \mathfrak{a})$,
the orbit $K_2\pi_1(x)$ is regular
if and only if $H$ is a regular point with respect to $(\tilde{\Sigma}, \Sigma, W)$.

All the symmetric triads with $\dim \mathfrak{a}=1$ are classified into the following four types (\cite{Ik1}):
\begin{center}
\begin{tabular}{|l|c|c|c|} 
\hline                            & $\Sigma^{+}$           & $W^{+}$               & $\tilde{\alpha}$ \\ \hline \hline
$\rm{III\mathchar`-B}_{1}$  & $\{\alpha \}$          & $\{\alpha \}$         & $\alpha$        \\ \hline
$\rm{I\mathchar`-BC}_{1}$   & $\{\alpha, 2\alpha \}$ & $\{\alpha \}$         & $\alpha$        \\ \hline
$\rm{II\mathchar`-BC}_{1}$  & $\{\alpha \}$          & $\{\alpha, 2\alpha \}$& $2\alpha$       \\ \hline
$\rm{III\mathchar`-BC}_{1}$ & $\{\alpha, 2\alpha \}$ & $\{\alpha, 2\alpha \}$& $2\alpha$       \\ 
\hline
\end{tabular}
\end{center}

\medskip\par
Let $\vartheta := \langle \tilde{\alpha}, H \rangle$ for $H \in \mathfrak{a}$.
Then, by (\ref{eq:cell}), $P_0 = \{ H \in \mathfrak{a} \mid 0 < \vartheta < \pi/2\}$ is a cell
in these types.
If $N_1$ is simply connected, then the orbit space of $K_2$-action on $N_1$ is identified with 
$\overline{P_0} = \{ H \in \mathfrak{a} \mid 0 \leq \vartheta \leq \pi/2\}$,
more precisely, each orbit meets $\pi_1(\exp \overline{P_0})$ at one point.
A point in the interior of the orbit space corresponds to a regular orbit,
and there exists a unique minimal (harmonic) orbit among regular orbits.
On the other hand, two endpoints of the orbit space correspond to singular orbits.
These singular orbits are minimal (harmonic), moreover these are weakly reflective (\cite{IST2009}).

In the following, we express the two equations $\| B_{H}\|^{2} = 1/2 $ and $\tau_{H} = 0$
in terms of $\vartheta$ for each type. 
For this purpose, here we should calculate $\langle \alpha, \alpha \rangle$.
\begin{align*}
\langle \alpha ,\alpha \rangle
&= -{\rm Killing}(\alpha ,\alpha) = -{\rm tr}(\mathrm{ad}(\alpha )^{2})\\
&= -\left\{ \sum_{\lambda \in \Sigma^{+}} \sum_{i=1}^{m(\lambda )} \langle \mathrm{ad}(\alpha)^{2} S_{\lambda ,i}, S_{\lambda ,i} \rangle \right. 
+  \sum_{\lambda \in \Sigma^{+}} \sum_{i=1}^{m(\lambda )} \langle \mathrm{ad}(\alpha)^{2} T_{\lambda ,i}, T_{\lambda ,i} \rangle  \\ 
&\left. +\sum_{\beta \in W^{+}} \sum_{j=1}^{n(\beta )} \langle \mathrm{ad}(\alpha )^{2} X_{\beta ,j}, X_{\beta ,j} \rangle 
 +\sum_{\beta \in W^{+}} \sum_{j=1}^{n(\beta )} \langle \mathrm{ad}(\alpha )^{2} Y_{\beta ,j}, Y_{\beta ,j} \rangle \right\}\\
&= \sum_{\lambda \in\Sigma^{+}}2m(\lambda) \langle \alpha ,\lambda \rangle^{2} +\sum_{\beta \in W^{+}}2n(\beta) \langle \alpha , \beta \rangle^{2}.  
\end{align*}
In the case of type $\rm{III\mathchar`-BC}_{1}$, i.e. $\Sigma^{+} = \{\alpha, 2\alpha \}$ and $W^{+} = \{\alpha, 2\alpha \}$, we can see that
\begin{align*}
\langle \alpha, \alpha \rangle
&= 2m(\alpha) \langle \alpha, \alpha \rangle^{2} + 2m(2\alpha) \langle \alpha, 2\alpha \rangle^{2}
+ 2n(\alpha) \langle \alpha, \alpha \rangle^{2} + 2n(2\alpha) \langle \alpha, 2\alpha \rangle^{2}\\
&= 2\langle \alpha, \alpha \rangle^{2}(m(\alpha)+4m(2\alpha ) + n(\alpha)+4n(2\alpha)).
\end{align*}
Therefore, we obtain
\begin{equation}\label{eq:norm_of_alpha}
\langle \alpha, \alpha \rangle = \frac{1}{2(m(\alpha ) + 4m(2\alpha ) + n(\alpha) + 4n(2\alpha))}.
\end{equation}
In the cases of other types,
we have $\langle \alpha, \alpha \rangle$ by letting $m(2\alpha) = 0$ (resp. $n(2\alpha) = 0$)
if $2\alpha \not\in \Sigma^+$ (resp. $2\alpha \not\in W^+$).

\subsection{Type $\rm{III\mathchar`-B}_{1}$}\label{2b1}
By (\ref{biheq}), the biharmonic condition $\| B_{H}\|^{2} = 1/2$ is equivalent to
$$
m(\alpha) + n(\alpha) = m(\alpha) (\cot \vartheta)^{2} + n(\alpha) (\tan \vartheta)^{2}
$$
for $H \in P_0$.
Thus we have
$$
\tan \vartheta = 1, \ \text{or} \ \sqrt{\frac{m(\alpha )}{n(\alpha )}}.
$$
On the other hand, by (\ref{harmeq}), the harmonic condition $\tau_{H}=0$ is equivalent to
$$
-m(\alpha) \cot \vartheta + n(\alpha) \tan \vartheta = 0.
$$
Thus we have
$$
\tan \vartheta = \sqrt{\frac{m(\alpha)}{n(\alpha)}}.
$$
By (2) of Theorem~\ref{Theorem 3.6},
the situation is divided into the following two cases:
\begin{enumerate}
\item When $m(\alpha) = n(\alpha)$, if an orbit $K_{2}\pi_{1}(x)$ is biharmonic, then it is harmonic.
\item When $m(\alpha) \neq n(\alpha)$, an orbit $K_{2}\pi_{1}(x)$ is proper biharmonic
if and only if $(\tan \vartheta)^2 = 1$ for $H \in P_0$.
In this case, a unique proper biharmonic orbit exists at the center of $P_0$,
namely $\vartheta = \pi/4$.
\end{enumerate}
 
\subsection{Type $\rm{I\mathchar`-BC}_{1}$}\label{1bc1}
We denote $m_{1} := m(\alpha)$, $m_{2} := m(2\alpha)$ and $n_1 := n(\alpha)$ for short.
Then, by (\ref{biheq}), the biharmonic condition $\| B_{H}\|^{2} = 1/2$ is equivalent to
$$
m_{1} + n_{1} + 4m_{2} = m_{1} (\cot \vartheta)^{2} + n_{1} (\tan \vartheta)^{2} + 4m_{2} (\cot 2\vartheta)^{2}.
$$
Thus, we have
$$
(\tan \vartheta)^2 = \frac{m_{1}+n_{1}+6m_{2} \pm \sqrt{(m_{1}+n_{1}+6m_{2})^2-4(n_{1}+m_{2})(m_{1}+m_{2})}}{2(n_{1}+m_{2})}.
$$
By (\ref{harmeq}), the harmonic condition $\tau_{H}=0$ is equivalent to
$$
-m_{1} \cot \vartheta + n_{1} \tan \vartheta - 4m_{2} \cot 2\vartheta = 0.
$$
Thus, we have
$$
(\tan \vartheta )^{2}= \frac{m_{1}+m_{2} }{n_{1}+m_{2}}.
$$
Since
\begin{align*}
0 &< \frac{m_{1}+n_{1}+6m_{2}-\sqrt{(m_{1}+n_{1}+6m_{2})^2-4(n_{1}+m_{2})(m_{1}+m_{2})}}{2(n_{1}+m_{2})}\\
&< \frac{m_{1}+m_{2}}{n_{1}+m_{2}}\\
&< \frac{m_{1}+n_{1}+6m_{2}+\sqrt{(m_{1}+n_{1}+6m_{2})^2-4(n_{1}+m_{2})(m_{1}+m_{2})}}{2(n_{1}+m_{2})},
\end{align*}
by (2) of Theorem~\ref{Theorem 3.6},
an orbit $K_{2}\pi_{1}(x)$ is proper biharmonic if and only if 
$$
(\tan \vartheta)^2 = \frac{m_{1}+n_{1}+6m_{2} \pm \sqrt{(m_{1}+n_{1}+6m_{2})^2-4(n_{1}+m_{2})(m_{1}+m_{2})}}{2(n_{1}+m_{2})}
$$
holds for $H \in P_0$.
Furthermore, a unique harmonic regular orbit exists between two proper biharmonic orbits in $P_0$.

\subsection{Type $\rm{II\mathchar`-BC}_{1}$}\label{2bc1}
By the definition of multiplicities, if $2\alpha \in W^+$, then $m(\alpha) = n(\alpha)$.
Hence we denote $m_1 := m(\alpha) = n(\alpha)$ and $n_2 := n(2\alpha)$. 
Then, by (\ref{biheq}), the biharmonic condition $\| B_{H}\|^{2} = 1/2$ is equivalent to
$$
2m_1+4n_2 = m_1 \big( (\cot (\vartheta/2))^{2} + (\tan (\vartheta/2))^{2} \big) + 4n_2 (\tan \vartheta)^{2}.
$$
Thus, we have
$$
(\tan \vartheta)^2 = \frac{n_2 \pm \sqrt{n_2^{2}-4n_2m_1}}{2n_2} = \frac{1}{2} \pm \sqrt{\frac{n_2-4m_1}{4n_2}}.
$$
By (\ref{harmeq}), the harmonic condition $\tau_{H}=0$ is equivalent to
$$
m_1 \big( -\cot (\vartheta/2) + \tan (\vartheta/2) \big) + 2n_2 \tan \vartheta = 0.
$$
Thus, we have
$$
(\tan \vartheta)^{2} =\frac{m_1}{n_2}.
$$
By (2) of Theorem~\ref{Theorem 3.6}, 
the situation is divided into the following three cases:
\begin{enumerate}
\item When $n_2<4m_1$, if $K_{2}\pi_{1}(x)$ is biharmonic, then it is harmonic.
\item When $n_2=4m_1$, an orbit $K_{2}\pi_{1}(x)$ is proper biharmonic if and only if $(\tan \vartheta)^2 = 1/2$
for $H \in P_0$. 
\item When $n_2>4m_1$, an orbit $K_{2}\pi_{1}(x)$ is proper biharmonic if and only if
$$
(\tan \vartheta)^2 = \frac{n_2 \pm \sqrt{n_2^{2}-4n_2m_1}}{2n_2}
$$
holds for $H \in P_0$, since
$$
0 < \frac{m_1}{n_2} < \frac{n_2-\sqrt{n_2^{2}-4n_2m_1}}{2n_2} < \frac{n_2+\sqrt{n_2^{2}-4n_2m_1}}{2n_2}.
$$
\end{enumerate}

\subsection{Type $\rm{III\mathchar`-BC}_{1}$}\label{3bc1}
By the definition of multiplicities, if $2\alpha \in W^+$, then $m(\alpha) = n(\alpha)$.
Hence we denote $m_{1} := m(\alpha) = n(\alpha), m_{2} := m(2\alpha)$ and $n_2 := n(2\alpha)$. 
Then, by (\ref{biheq}), the biharmonic condition $\| B_{H}\|^{2} = 1/2$ is equivalent to
$$
2m_{1}+4m_{2}+4n_2 = m_{1} \big( (\cot (\vartheta/2))^{2}+(\tan (\vartheta/2))^{2} \big)+4m_{2}(\cot \vartheta)^{2}+4n_2(\tan \vartheta)^{2}.
$$
Thus, we have
\begin{align*}
(\tan \vartheta)^{2}
&= \frac{m_{2}+n_2 \pm \sqrt{(m_{2}+n_2)^2-4n_2(m_{1}+m_{2})}}{2n_2} \\
&= \frac{m_{2}+n_2 \pm \sqrt{(m_{2}-n_2)^2-4n_2m_{1}}}{2n_2}.
\end{align*}
By (\ref{harmeq}), the harmonic condition $\tau_{H}=0$ is equivalent to
$$
m_{1} \big(\tan (\vartheta/2) - \cot (\vartheta/2) \big) - 2 m_{2} \cot \vartheta + 2n_2 \tan \vartheta = 0.
$$
Thus, we have
$$
(\tan \vartheta)^2 = \frac{m_{1}+m_{2}}{n_2}.
$$
By (2) of Theorem~\ref{Theorem 3.6}, we obtain the following results:
\begin{enumerate}
\item When $(m_{2}-n_2)^2 - 4n_2m_{1}<0$, if $K_{2}\pi_{1}(x)$ is biharmonic, then it is harmonic.
\item When $(m_{2}-n_2)^2 - 4n_2m_{1}=0$, an orbit $K_{2}\pi_{1}(x)$ is proper biharmonic if and only if
$(\tan \vartheta)^2 = (m_{2}+n_2)/2n_2$ for $H \in P_0$.
\item When $(m_{2}-n_2)^2 - 4n_2m_{1}>0$, an orbit $K_{2}\pi_{1}(x)$ is proper biharmonic if and only if 
$$
(\tan \vartheta)^{2} = \frac{m_{2} + n_2 \pm \sqrt{(m_{2}-n_2)^2 -4n_2m_{1}}}{2n_2}
$$
for $H \in P_0$.
\end{enumerate}

For the proof of (2), we will show that 
$$
\frac{m_{1}+m_{2}}{n_2} \neq \frac{m_{2}+n_2}{2n_2}.
$$
If $(m_{1}+m_{2})/n_2 =(m_{2}+n_2)/(2n_2)$, then $2m_{1}+m_{2}-n_2=0$.
Hence $(m_{2}-n_2)^{2}-4n_2m_{1}=-4m_{1}(m_{1}+m_{2})<0$, 
which is a contradiction.

For the proof of (3), we will show that 
$$
\frac{m_{1}+m_{2}}{n_2} \neq \frac{m_{2}+n_2\pm \sqrt{(m_{2}-n_2)^{2}-4n_2m_{1}}}{2n_2}.
$$
If the equality holds, then we have $(2m_{1}+m_{2}-n_2)^{2}=(m_{2}-n_2)^{2}-4n_2m_{1}$.
Hence $4m_{1}(m_{1}+m_{2}) =0$, which is a contradiction.

In fact, in the cases of type $\rm{III\mathchar`-BC}_{1}$, a compact symmetric triad which is not (1) 
is only $(E_{6}, \mathrm{SO}(10)\cdot\mathrm{U}(1), F_{4})$ in the list below.
In this case, 
\begin{align*}
\frac{m_{1}+m_{2}}{n_2}
&< \frac{m_{2}+n_2-\sqrt{(m_{2}-n_2)^{2}-4n_2m_{1}}}{2n_2} \\
&< \frac{m_{2}+n_2+ \sqrt{(m_{2}-n_2)^{2}-4n_2m_{1}}}{2n_2}
\end{align*}
holds.

\vskip0.6cm\par
Let $b>0$, $c>1$ and $q>1$.
Each commutative compact symmetric triad $(G, K_{1}, K_{2})$
where $G$ is simple, $\theta_1 \not\sim \theta_2$ and $\dim \mathfrak{a}=1$
is one of the following (see \cite{Ik2}):

\smallskip
\noindent
{\bf Type $\rm{III\mathchar`-B}_{1}$}\\
\begin{tabular}{|c|c|} 
\hline
$(G,K_{1},K_{2})$                                                                   &$(m(\alpha),n(\alpha))$  \\ \hline \hline
$(\mathrm{SO}(1+b+c), \mathrm{SO}(1+b) \times \mathrm{SO}(c), \mathrm{SO}(b+c))$    &$(c-1,b)$\\ \hline
$(\mathrm{SU}(4), \mathrm{Sp}(2), \mathrm{SO}(4))$                                  &$(2,2)$  \\ \hline
$(\mathrm{SU}(4), \mathrm{S}(\mathrm{U}(2) \times \mathrm{U}(2)), \mathrm{Sp}(2))$ &$(3,1)$  \\ \hline
$(\mathrm{Sp}(2), \mathrm{U}(2), \mathrm{Sp}(1) \times \mathrm{Sp}(1))$             &$(1,2)$  \\ \hline
\end{tabular}

\smallskip
\noindent
{\bf Type $\rm{I\mathchar`-BC}_{1}$}\\
\begin{tabular}{|c|c|} 
\hline
$(G,K_{1},K_{2})$                                                                           &$(m(\alpha),m(2\alpha),n(\alpha))$  \\ \hline \hline
$(\mathrm{SO}(2+2q), \mathrm{SO}(2) \times \mathrm{SO}(2q), \mathrm{U}(1+q))$               &$(2(q-1),1,2(q-1))$\\ \hline
$(\mathrm{SU}(1+b+c), \mathrm{S}(\mathrm{U}(1+b) \times \mathrm{U}(c)), \mathrm{S}(\mathrm{U}(1) \times \mathrm{U}(b+c))$&$(2(c-1),1,2b)$\\ \hline
$(\mathrm{Sp}(1+b+c), \mathrm{Sp}(1+b) \times \mathrm{Sp}(c), \mathrm{Sp}(1) \times \mathrm{Sp}(b+c))$&$(4(c-1),3,4b)$\\ \hline
$(\mathrm{SO}(8), \mathrm{U}(4), \mathrm{U}(4)')$                                             &$(4,1,1)$ \\ \hline
\end{tabular}

\smallskip
\noindent
{\bf Type $\rm{II\mathchar`-BC}_{1}$}\\
\begin{tabular}{|c|c|} 
\hline
$(G,K_{1},K_{2})$                                                                           &$(m(\alpha),n(\alpha),n(2\alpha))$  \\ \hline \hline
$(\mathrm{SO}(6), \mathrm{U}(3), \mathrm{SO}(3) \times \mathrm{SO}(3))$                         &$(2,2,1)$ \\ \hline
$(\mathrm{SU}(1+q), \mathrm{SO}(1+q), \mathrm{S}(\mathrm{U}(1) \times \mathrm{U}(q)))$   &$(q-1,q-1,1)$\\ \hline
\end{tabular}

\smallskip
\noindent
{\bf Type $\rm{III\mathchar`-BC}_{1}$}\\
\begin{tabular}{|c|c|} 
\hline
$(G,K_{1},K_{2})$                                                                           &$(m(\alpha),m(2\alpha),n(\alpha),n(2\alpha))$  \\ \hline \hline
$(\mathrm{SU}(2+2q), \mathrm{S}(\mathrm{U}(2) \times \mathrm{U}(2q)), \mathrm{Sp}(1+q))$     &$(4(q-1),3,4(q-1),1)$\\ \hline
$(\mathrm{Sp}(1+q), \mathrm{U}(1+q), \mathrm{Sp}(1) \times \mathrm{Sp}(q)) $          &$(2(q-1),1,2(q-1),2)$\\ \hline
$(\mathrm{E}_{6}, \mathrm{SU}(6) \cdot \mathrm{SU}(2), \mathrm{F}_{4})$                         &$(8,3,8,5)$ \\ \hline
$(\mathrm{E}_{6}, \mathrm{SO}(10) \cdot \mathrm{U}(1), \mathrm{F}_{4})$                         &$(8,7,8,1)$ \\ \hline
$(\mathrm{F}_{4}, \mathrm{Sp}(3) \cdot \mathrm{Sp}(1), \mathrm{Spin}(9))$                         &$(4,3,4,4)$\\ \hline
\end{tabular}

\medskip
Here, we define 
$\mathrm{U}(4)'=\{g\in \mathrm{SO}(8) \mid JgJ^{-1}=g \}$ 
where
$$
J=\left[
\begin{array}{cc|cc}
& &I_{3} &\\
& & &-1\\ \hline
-I_{3}& & &\\
& 1& &
\end{array}
\right]
$$ 
and $I_l$ denotes the identity matrix of $l \times l$.

\vskip0.6cm\par
 
\section{Main result and examples} \label{sect:results}
{\bf 6.1.} 
Summing up the previous sections, 
we classify all the biharmonic hypersurfaces in irreducible compact symmetric spaces
which are orbits of commutative Hermann actions,
namely we obtain the following theorem.
\begin{theorem}\label{thm:list_of_biharmonic_orbits}
Let $(G,K_1,K_2)$ be a commutative compact symmetric triad where $G$ is simple,
and suppose that $K_2$-action on $N_1=G/K_1$ is cohomogeneity one
(hence $K_1$-action on $N_2=G/K_2$ is also cohomogeneity one).
Then all the proper biharmonic hypersurfaces which are regular orbits of $K_2$-action (resp. $K_1$-action)
in the compact symmetric space $N_{1}$ (resp. $N_{2}$)
are classified into the following lists:
\begin{enumerate}
\item When $(G,K_1,K_2)$ is one of the following cases,
there exists a unique proper biharmonic hypersurface which is a regular orbit
of $K_{2}$-action on $N_{1}$
(resp. $K_{1}$-action on $N_{2}$).
\begin{enumerate}
\item[(1-1) ] $(\mathrm{SO}(1+b+c),\ \mathrm{SO}(1+b) \times \mathrm{SO}(c),\ \mathrm{SO}(b+c)) \quad (b>0,\ c>1,\ c-1\neq b)$
\item[(1-2) ] $(\mathrm{SU}(4),\ \mathrm{S}(\mathrm{U}(2) \times \mathrm{U}(2)),\ \mathrm{Sp}(2))$
\item[(1-3) ] $(\mathrm{Sp}(2),\ \mathrm{U}(2),\ \mathrm{Sp}(1) \times \mathrm{Sp}(1))$ 
\end{enumerate}

\item When $(G,K_1,K_2)$ is one of the following cases,
there exist exactly two distinct proper biharmonic hypersurfaces which are regular orbits of
of $K_{2}$-action on $N_{1}$
(resp. $K_{1}$-action on $N_{2}$).
\begin{enumerate}
\item[(2-1) ] $(\mathrm{SO}(2+2q),\ \mathrm{SO}(2) \times \mathrm{SO}(2q),\ \mathrm{U}(1+q)) \quad (q>1)$
\item[(2-2) ] $(\mathrm{SU}(1+b+c),\ \mathrm{S}(\mathrm{U}(1+b) \times \mathrm{U}(c)),\ \mathrm{S}(\mathrm{U}(1) \times \mathrm{U}(b+c)) \quad (b \geq 0,\ c>1)$
\item[(2-3) ] $(\mathrm{Sp}(1+b+c),\ \mathrm{Sp}(1+b) \times \mathrm{Sp}(c),\ \mathrm{Sp}(1) \times \mathrm{Sp}(b+c)) \quad (b \geq 0,\ c>1)$
\item[(2-4) ] $(\mathrm{SO}(8),\ \mathrm{U}(4),\ \mathrm{U}(4)')$
\item[(2-5) ] $(\mathrm{E}_{6},\ \mathrm{SO}(10) \cdot \mathrm{U}(1),\ \mathrm{F}_{4})$
\item[(2-6) ] $(\mathrm{SO}(1+q),\ \mathrm{SO}(q),\ \mathrm{SO}(q)) \quad (q>1)$
\item[(2-7) ] $(\mathrm{F}_4,\ \mathrm{Spin}(9),\ \mathrm{Spin}(9))$
\end{enumerate}

\item When $(G,K_1,K_2)$ is one of the following cases,
any biharmonic regular orbit of $K_{2}$-action on $N_{1}$
(resp. $K_{1}$-action on $N_{2}$) is harmonic.
\begin{enumerate}
\item[(3-1) ] $(\mathrm{SO}(2c),\ \mathrm{SO}(c) \times \mathrm{SO}(c),\ \mathrm{SO}(2c-1)) \quad (c>1)$
\item[(3-2) ] $(\mathrm{SU}(4),\ \mathrm{Sp}(2),\ \mathrm{SO}(4))$
\item[(3-3) ] $(\mathrm{SO}(6),\ \mathrm{U}(3),\ \mathrm{SO}(3) \times \mathrm{SO}(3))$
\item[(3-4) ] $(\mathrm{SU}(1+q),\ \mathrm{SO}(1+q),\ \mathrm{S}(\mathrm{U}(1) \times \mathrm{U}(q))) \quad (q>1)$
\item[(3-5) ] $(\mathrm{SU}(2+2q),\ \mathrm{S}(\mathrm{U}(2) \times \mathrm{U}(2q)),\ \mathrm{Sp}(1+q)) \quad (q>1)$
\item[(3-6) ] $(\mathrm{Sp}(1+q),\ \mathrm{U}(1+q),\ \mathrm{Sp}(1) \times \mathrm{Sp}(q)) \quad (q>1)$
\item[(3-7) ] $(\mathrm{E}_{6},\ \mathrm{SU}(6) \cdot \mathrm{SU}(2),\ \mathrm{F}_{4})$
\item[(3-8) ] $(\mathrm{F}_{4},\ \mathrm{Sp}(3) \cdot \mathrm{Sp}(1),\ \mathrm{Spin}(9))$
\end{enumerate}
\end{enumerate}
\end{theorem}

\begin{remark}
In Theorem~\ref{thm:list_of_biharmonic_orbits},
we determined all the biharmonic hypersurfaces in irreducible compact symmetric spaces
which are orbits of commutative Hermann actions.
\begin{enumerate}
\item In the previous section we assumed $\theta_1 \not\sim \theta_2$.
If $\theta_1 \sim \theta_2$, then the action of $K_2$ on $N_1$ is orbit equivalent
to the isotropy action of $K_1$ on $N_1$.
We will discuss these cases in Section 6.3.
\item The commutative condition $\theta_1 \theta_2 = \theta_2\theta_1$ is essential
for our discussion.
Indeed, there exist some Hermann actions where $\theta_1 \theta_2 \neq \theta_2\theta_1$.
Moreover there exist some hyperpolar actions of cohomogeneity one on irreducible compact symmetric spaces
which are not Hermann actions (cf. \cite{Kollross}).
\end{enumerate}
\end{remark}

\vskip0.6cm\par
{\bf 6.2.} 
We shall explain details of the cases (1-1), (2-2) and (3-1) in Theorem~\ref{thm:list_of_biharmonic_orbits},
and give new examples of proper biharmonic orbits.
By Proposition~\ref{prop5.1}, a proper biharmonic orbit $K_{2}\pi_{1}(x)$ in $N_{1}$ corresponds to a proper biharmonic orbit 
$K_1\pi_{2}(x)$ in $N_{2}$.
In particular, we can obtain new examples of proper biharmonic orbits
corresponding to some known examples.

We consider the isotropy subgroups of orbits of Hermann actions.
For $x= \exp H \ (H \in \mathfrak{a})$, we define the isotropy subgroups
$$
(K_{2})_{\pi_{1}(x)}=\{k \in K_{2} \mid k\pi_{1}(x)=\pi_{1}(x)\}, 
$$
$$
(K_{1})_{\pi_{2}(x)}=\{k \in K_{1} \mid k\pi_{2}(x)=\pi_{2}(x)\}.
$$
Then we can show that $(K_{2})_{\pi_{1}(x)} \cong (K_{1})_{\pi_{2}(x)}$ by an inner automorphism of $G$.
The orbit $K_2\pi_1(x)$ (resp. $K_1\pi_2(x)$) is diffeomorphic to the homogeneous space $K_2/((K_{2})_{\pi_{1}(x)})$ (resp. $K_1/((K_{1})_{\pi_{2}(x)})$).
If $K_{2}\pi_{1}(x)$ is a regular orbit, 
then $K_{1}\pi_{2}(x)$ is also a regular orbit, and we have
$\mathrm{Lie}((K_{2})_{\pi_{1}(x)}) = \mathrm{Lie}((K_{1})_{\pi_{2}(x)})= \mathfrak{k}_{0}.$ 

\vskip0.6cm\par
\noindent
{\bf Example 1.}\quad
Let $(G,K_{1},K_{2})=(\mathrm{SO}(1+b+c), \mathrm{SO}(1+b) \times \mathrm{SO}(c), \mathrm{SO}(b+c))\ (b>0,\ c>1).$
This is the case of (3-1) when $c-1=b$,
otherwise the case of (1-1) in Theorem~\ref{thm:list_of_biharmonic_orbits}.
In this case, the involutions $\theta_1$ and $\theta_2$ are given by
$$
\theta_{1}(k)=I'_{1+b}k I'_{1+b},\quad \theta_{2}(k)=I'_{1}k I'_{1} \quad (k \in G), 
$$
where 
$$
I'_{l}=\left[
\begin{array}{cc}
-I_{l}&0\\
0&I_{1+b+c-l}
\end{array}
\right]
\ \ (1\leq l\leq b+c).
$$
Then, we have the canonical decompositions $\mathfrak{g} = \mathfrak{k}_1 \oplus \mathfrak{m}_1 = \mathfrak{k}_2 \oplus \mathfrak{m}_2$ as
\begin{align*}
\mathfrak{k}_{1} &=
\left\{
\left[
\begin{array}{cc}
X&0\\
0&Y
\end{array}
\right]
\ \bigg| \ 
{X \in \mathfrak{so}(1+b) \atop  Y \in \mathfrak{so}(c)}
\right\},&
\mathfrak{m}_{1} &=
\left\{
\left[
\begin{array}{cc}
0&X\\
-{}^{t}X&0
\end{array}
\right]
\ \bigg| \ 
{X \in \mathrm{M}_{1+b, c}(\mathbb{R})}
\right\}, \\
\mathfrak{k}_{2} &=
\left\{
\left[
\begin{array}{cc}
0&0\\
0&X
\end{array}
\right]
\ \bigg| \ 
X \in \mathfrak{so}(b+c) 
\right\},&
\mathfrak{m}_{2} &=
\left\{ 
\left[
\begin{array}{cc}
0&X\\
-{}^{t}X&0
\end{array}
\right]
\ \bigg| \ 
X \in \mathrm{M}_{1, b+c}(\mathbb{R})
\right\}.
\end{align*}
Thus, we have 
\begin{align*}
\mathfrak{k}_{1} \cap \mathfrak{k}_{2} &=
\left\{
\left[
\begin{array}{ccc}
0&0&0\\
0&X&0\\
0&0&Y
\end{array}
\right]
\ \Bigg| \ 
{X \in \mathfrak{so}(b) \atop Y \in \mathfrak{so}(c)}
\right\}, \\
\mathfrak{m}_{1} \cap \mathfrak{m}_{2} &=
\left\{ 
\left[
\begin{array}{ccc}
0&0&X\\
0&0&0\\
-{}^{t}X&0&0
\end{array}
\right]
\ \Bigg| \ 
X \in \mathrm{M}_{1,c}(\mathbb{R})
\right\}, \\
\mathfrak{k}_{1} \cap \mathfrak{m}_{2} &=
\left\{
\left[
\begin{array}{ccc}
0&X&0\\
-{}^{t}X&0&0\\
0&0&0
\end{array}
\right]
\ \Bigg| \ 
X \in \mathrm{M}_{1,b}(\mathbb{R})
\right\}, \\
\mathfrak{m}_{1} \cap \mathfrak{k}_{2} &=
\left\{
\left[
\begin{array}{ccc}
0&0&0\\
0&0&X\\
0&-{}^{t}X&0
\end{array}
\right]
\ \Bigg| \ 
X \in \mathrm{M}_{b,c}(\mathbb{R})
\right\}.
\end{align*}
We take a maximal abelian subspace $\mathfrak{a}$ in $\mathfrak{m}_{1} \cap \mathfrak{m}_{2}$ as
$$
\mathfrak{a}=
\left\{
H(\vartheta) = 
\left[
\begin{array}{ccc}
0&0&X\\
0&0&0\\
-{}^{t}X&0&0
\end{array}
\right]
\ \Bigg| \ 
{X = [0,\ldots,0,\vartheta] \atop
\vartheta \in \mathbb{R}}
\right\}.
$$
Then we have 
\begin{align*}
\mathfrak{k}_{0}&=
\left\{ 
\left.
\left[
\begin{array}{cccc}
0&0&0&0\\
0&X&0&0\\
0&0&Y&0\\
0&0&0&0
\end{array}
\right]
\ \right|\ 
{X \in \mathfrak{so}(b) \atop Y \in \mathfrak{so}(c-1)}
\right\}
,\\
V(\mathfrak{k}_{1}\cap \mathfrak{m}_{2})&=\{ 0\}
,\\
V(\mathfrak{m}_{1}\cap \mathfrak{k}_{2})&=
\left\{ 
\left.
\left[
\begin{array}{cccc}
0&0&0&0\\
0&0&X&0\\
0&-{}^{t}X&0&0\\
0&0&0&0
\end{array}
\right]
\ \right|\ 
X \in \mathrm{M}_{b,c-1}(\mathbb{R})
\right\}.
\end{align*}

Let $E_{i}^{j}$ be a matrix whose $(i,j)$-entry is one and all the other entries are zero.
We define $A_{i}^{j}:=E_{i}^{j}- E_{j}^{i}$. 
Then, we can see
\begin{align*}
[H(\vartheta), A_{1}^{j}] &= -\vartheta A_{1+b+c}^{j} &(2\leq j\leq b+c),\\
[H(\vartheta), A_{1+b+c}^{j}] &= \vartheta A_{1}^{j} &(2\leq j\leq b+c).
\end{align*}
We define a vector $\alpha \in \mathfrak{a}$ by $\langle H(\vartheta), \alpha \rangle = \vartheta \ (\vartheta \in \mathbb{R})$.
Then
\begin{align*}
\mathfrak{k}_{\alpha} &= \mathrm{Span} \{A_{1+b+c}^{2+b},\ldots , A_{1+b+c}^{b+c}\},\\
\mathfrak{m}_{\alpha} &= \mathrm{Span} \{A_{1}^{2+b},\ldots , A_{1}^{b+c}\},\\
V_{\alpha}^{\perp}( \mathfrak{k}_{1}\cap \mathfrak{m}_{2}) &= \mathrm{Span} \{A_{1}^{2},\ldots , A_{1}^{1+b}\},\\
V_{\alpha}^{\perp}( \mathfrak{m}_{1}\cap \mathfrak{k}_{2}) &= \mathrm{Span} \{A_{1+b+c}^{2},\ldots , A_{1+b+c}^{1+b}\}.
\end{align*}
Hence, in this case, we have
$$
\Sigma^{+} = \{\alpha\},\quad
W^{+} = \{\alpha\},\quad
m(\alpha) = c-1,\quad
n(\alpha) = b.
$$ 

Let $x_{0}=\exp (H(\pi/4))$.
By the computation in Section~\ref{2b1}, we can see that 
$K_{2}\pi_{1}(x_{0})$ and $K_1\pi_{2}(x_{0})$ are biharmonic hypersurfaces in $N_1$ and $N_2$, respectively.
These orbits exist at the center of the orbit space $\overline{P_0} = \{H(\vartheta) \mid 0 \leq \vartheta \leq \pi/2\}$.
When $c-1=b$, these orbits are harmonic.
When $c-1\neq b$, these are not harmonic, hence proper biharmonic.
The orbit $K_2\pi_{1}(x_0)$ is the Clifford hypersurface
$S^{b}(1/{\sqrt 2})\times S^{c-1}(1/\sqrt{2}) \cong (\mathrm{SO}(1+b) \times \mathrm{SO}(c)) / (\mathrm{SO}(b) \times \mathrm{SO}(c-1))$
embedded in the sphere $S^{b+c}(1) \cong \mathrm{SO}(1+b+c)/\mathrm{SO}(b+c) = N_2$ (\cite{J}).
On the other hand, the orbit $K_{2}\pi_{1}(x_{0})$ is diffeomorphic to $\mathrm{SO}(b+c)/(\mathrm{SO}(b) \times \mathrm{SO}(c-1))$,
i.e. the universal covering of a real flag manifold, and embedded in the oriented real Grassmannian manifold 
$\widetilde{G_{1+b}}(\mathbb{R}^{1+b+c}) \cong \mathrm{SO}(1+b+c)/(\mathrm{SO}(1+b) \times \mathrm{SO}(c))=N_{1}$
as the tube of radius $\pi/4$ over the totally geodesic sub-Grassmannian $\widetilde{G_{b}}(\mathbb{R}^{b+c})$.
The orbit $K_{2}\pi_{1}(x_{0})$ in $N_1$ gives a new example of a proper biharmonic hypersurface
in the oriented real Grassmannian manifold.

\vskip0.6cm\par
\noindent
{\bf Example 2.}\quad
Let $(G,K_{1},K_{2})=(\mathrm{SU}(1+b+c), 
\mathrm{S}(\mathrm{U}(1+b) \times \mathrm{U}(c)),
\mathrm{S}(\mathrm{U}(1) \times \mathrm{U}(b+c)))\
(b>0,\ c>1).$
This is the case of (2-2) except for $b=0$ in Theorem~\ref{thm:list_of_biharmonic_orbits}.
In this case, the involutions $\theta_1$ and $\theta_2$ are given by
$$
\theta_{1}(k) = I'_{1+b}k I'_{1+b},\quad \theta_{2}(k) = I'_{1}k I'_{1} \quad (k \in G).
$$
Analogous to the previous example, in this case, we have
$$
\Sigma^{+}=\{\alpha, 2\alpha \},\ W^{+}=\{\alpha\},\ m(\alpha)=2(c-1),\ m(2\alpha)=1,\ n(\alpha)=2b.
$$ 
Therefore, the symmetric triad $(\tilde{\Sigma}, \Sigma, W)$ is of type $\rm{I\mathchar`-BC}_{1}$. 
By the computation in Section~\ref{1bc1},
we have two distinct proper biharmonic hypersurfaces in $N_{1}$, and also in $N_{2}$.
More precisely, 
let $x_\pm = \exp(H(\vartheta_\pm))$ where $\vartheta_\pm$ is a solution of the equation
\begin{align*}
(\tan \vartheta) ^2&=\frac{m_{1}+n_{1} +6m_{2} \pm \sqrt{(m_{1}+n_{1} +6m_{2})^2- 4(n_{1} +m_{2})(m_{1} +m_{2}) } }{2(n_{1} + m_{2} )}\\
&=\frac{(c-1)+b +3 \pm \sqrt{((c-1)+b +3)^2- (2b +1)(2(c-1) +1) } }{2b + 1}.
\end{align*}
Then $K_{2}\pi_{1}(x_\pm)$ and $K_1\pi_{2}(x_\pm)$ are proper biharmonic hypersurfaces in $N_1$ and $N_2$, respectively.
The orbit $K_1\pi_2(x_\pm) \cong \mathrm{S}(\mathrm{U}(1+b) \times \mathrm{U}(c)) / \mathrm{S}(\mathrm{U}(b) \times \mathrm{U}(c-1) \times \mathrm{U}(1))$
is the tube of radius $\vartheta_\pm$ over the totally geodesic $\mathbb{C}P^{b}$
in the complex projective space $\mathbb{C}P^{b+c} \cong \mathrm{SU}(1+b+c) / \mathrm{S}(\mathrm{U}(1) \times \mathrm{U}(b+c)) = N_2$ (see Theorem 5 in \cite{IIU}).
On the other hand,
the orbit $K_2\pi_{1}(x_\pm) \cong \mathrm{S}(\mathrm{U}(1) \times \mathrm{U}(b+c)) / \mathrm{S}(\mathrm{U}(b) \times \mathrm{U}(c-1) \times \mathrm{U}(1))$
is the tube of radius $\vartheta_\pm$ over the totally geodesic sub-Grassmannian $G_{b}(\mathbb{C}^{b+c})$
in the complex Grassmannian manifold $G_{1+b}(\mathbb{C}^{1+b+c}) \cong \mathrm{SU}(1+b+c) / \mathrm{S}(\mathrm{U}(1+b) \times \mathrm{U}(c)) = N_1$.
The orbit $K_{2}\pi_{1}(x_\pm)$ in $N_1$ gives a new example of a proper biharmonic hypersurface
in the complex Grassmannian manifold.

\vskip0.6cm\par
{\bf 6.3.} 
In the above argument, we supposed that $\theta_{1}\not\sim \theta_{2}$
in order to use the classification of commutative compact symmetric triads.
However, we can apply our method to the cases of $\theta_{1} \sim \theta_{2}$.
When $\theta_{1} \sim \theta_{2}$, a Hermann action is orbit equivalent to the isotropy action of a compact symmetric space (see  \cite{Ik1}). 
Hence, it is sufficient to discuss the cases of isotropy actions, that is, $\theta_{1}=\theta_{2}$.
When $\theta_{1}=\theta_{2}$, we have $W=\emptyset$, since $\mathfrak{k}_{1}\cap \mathfrak{m}_{2}=\mathfrak{m}_{1}\cap \mathfrak{k}_{2}=\{0\}$.
Thus we have $\tilde{\Sigma}=\Sigma$.
Moreover, $\Sigma$ is the root system of the compact symmetric space $N_{1}$ with respect to $\mathfrak{a}$.
Since we consider the cases of $\dim \mathfrak{a}=1$, 
the rank of $N_{1}$ equals to one.
All the simply connected, rank one symmetric spaces of compact type are classified as follows:
$$
S^{q},\ \mathbb{C}P^{q},\ \mathbb{H}P^{q},\ \mathbb{O}P^{2} \quad (q\geq 2).
$$
The isotropy actions of these symmetric spaces correspond to the cases
(2-6), (2-2) with $b=0$, (2-3) with $b=0$, and (2-9) in Theorem~\ref{thm:list_of_biharmonic_orbits},
respectively.
Except for the case of $\mathbb{O}P^{2}$, homogeneous biharmonic hypersurfaces in compact, rank one symmetric spaces were classified (\cite{IIU2}, \cite{IIU}).
Therefore, we consider the octonionic projective plane $\mathbb{O}P^{2}\cong F_{4}/\mathrm{Spin}(9)$.

Let $(G, K_{1}, K_{2})=(F_{4}, \mathrm{Spin}(9),\mathrm{Spin}(9))$ with $\theta_{1}=\theta_{2}$.
This is the case of (2-9) in Theorem~\ref{thm:list_of_biharmonic_orbits}.
Since $K_{1}=K_{2}$, we denote 
$$
\mathfrak{k}:=\mathfrak{k}_{1}=\mathfrak{k}_{2},\quad
\mathfrak{m}:=\mathfrak{m}_{1}=\mathfrak{m}_{2}.
$$
We define an $\mathrm{Ad}(G)$-invariant inner product on ${\mathfrak g}$
by $\langle \cdot , \cdot \rangle= -\mathrm{Killing}(\cdot ,\cdot)$.
Fix a maximal abelian subspace ${\mathfrak a}$ in ${\mathfrak m}$.
Then we have $\Sigma^{+}=\{\alpha, 2\alpha \}$ and $m(\alpha)=8, \ m(2\alpha)=7$ 
(\cite{He}, Page 534).
By letting $n(\alpha) = n(2\alpha) = 0$ in (\ref{eq:norm_of_alpha}) since $W^+=\emptyset$,
we can see that $\langle \alpha, \alpha \rangle = 1/\{2(8+4\cdot 7)\}$.
Let $x = \exp H$ for $H \in \mathfrak{a}$.
By Theorem~\ref{thm:norm_of_2nd_fundamental_form1}, 
we have the following:
\begin{align*}
\|B_{H}\|^{2}&=8(\cot \langle \alpha, H\rangle )^{2} \langle \alpha, \alpha\rangle +
7(\cot \langle 2\alpha, H\rangle )^{2} \langle 2\alpha, 2\alpha \rangle \\
&= \frac{1}{18}\{ 2(\cot \langle \alpha, H\rangle )^{2} +7(\cot \langle 2\alpha, H\rangle )^{2} \} ,\\
dL_{x}^{-1}(\tau_{H})&= -(8\cot \langle \alpha, H\rangle \alpha + 7\cot \langle 2\alpha, H\rangle 2\alpha )\\
&=-2(4\cot \langle \alpha, H\rangle + 7\cot \langle 2\alpha, H\rangle ) \alpha .
\end{align*} 
Then, the biharmonic condition $\|B_{H}\|^{2}=1/2$ is equivalent to 
$$
9= 2(\cot \langle \alpha, H\rangle )^{2} +7(\cot \langle 2\alpha, H\rangle )^{2}.
$$
Thus we have 
$$
(\cot \langle \alpha, H\rangle )^{2}=\frac{25 \pm 2\sqrt{130}}{15}.
$$
The harmonic condition $\tau_{H}=0$ is equivalent to 
$$
4\cot \langle \alpha, H\rangle + 7\cot \langle 2\alpha, H\rangle =0.
$$
Thus we have 
$$
(\cot \langle \alpha, H\rangle )^{2}=\frac{7}{15}.
$$
Since 
$$
0<\frac{25 - 2\sqrt{130}}{15}< \frac{7}{15}<\frac{25+ 2\sqrt{130}}{15},
$$
by (2) of Theorem~\ref{Theorem 3.6}, 
an orbit $K_{2}\pi_{1}(x)$ is proper biharmonic if and only if 
$$
(\cot \langle \alpha, H\rangle )^{2}=\frac{25 \pm 2\sqrt{130}}{15}
$$
holds for $H \in \mathfrak{a}$ with $0< \langle \alpha, H\rangle <\pi/2$.
Furthermore, a unique harmonic regular orbit exists between two proper biharmonic orbits in $\{H \in \mathfrak{a} \mid 0< \langle \alpha, H\rangle <\pi/2 \}$.
These regular orbits are diffeomorphic to $S^{15}$ embedded in $\mathbb{O}P^2$.

\subsection*{Acknowledgement}
The authors would like to thank Professor Osamu Ikawa for his useful comments.

\bibliographystyle{amsplain}

\end{document}